# Data-driven wavelet-Fisz methodology for nonparametric function estimation

## Piotr Fryzlewicz

*Department of Mathematics*
*University of Bristol*
*University Walk*
*Bristol BS8 1TW*
*UK*
*e-mail:* p.z.fryzlewicz@bris.ac.uk

**Abstract:** We propose a wavelet-based technique for the nonparametric estimation of functions contaminated with noise whose mean and variance are linked via a possibly unknown variance function. Our method, termed the data-driven wavelet-Fisz technique, consists of estimating the variance function via a Nadaraya-Watson estimator, and then performing a wavelet thresholding procedure which uses the estimated variance function and local means of the data to set the thresholds at a suitable level.

We demonstrate the mean-square near-optimality of our wavelet estimator over the usual range of Besov classes. To achieve this, we establish an exponential inequality for the Nadaraya-Watson variance function estimator.

We discuss various implementation issues concerning our wavelet estimator, and demonstrate its good practical performance. We also show how it leads to a new wavelet-domain data-driven variance-stabilising transform. Our estimator can be applied to a variety of problems, including the estimation of volatilities, spectral densities and Poisson intensities, as well as to a range of problems in which the distribution of the noise is unknown.

Received November 2007.

## 1. Introduction

A paradigmatic problem in nonparametric regression is the estimation of a one-dimensional function $\alpha : [0, 1] \mapsto \mathbb{R}$ from noisy observations $X_t$ taken on an equispaced grid:

$$X_t = \alpha(t/n) + \varepsilon_t, \quad t = 1, \ldots, n, \tag{1}$$

where the $\varepsilon_t$'s are random variables with $\mathbb{E}(\varepsilon_t) = 0$. Various subclasses of the problem can be identified, depending on the smoothness properties of $\alpha$ and the joint distribution of $(\varepsilon_t)_{t=1}^n$.

Since the seminal work of Donoho and Johnstone (1994), estimation techniques based on non-linear wavelet shrinkage have become a commonly used tool in nonparametric function estimation, extensively studied in the statistical literature. Many of them combine excellent finite-sample performance, linear computational complexity, and optimal (or near-optimal) asymptotic Mean-Square Error behaviour over a variety of function smoothness classes, including functions of a low degree of regularity. This often puts them at an advantage compared to linear estimation techniques (such as kernel smoothing) in setups where the function $\alpha$ has discontinuities or exhibits an otherwise irregular behaviour. A comprehensive overview of wavelet methods in statistics can be found, for example, in Vidakovic (1999).

The main idea underlying most wavelet techniques is that upon transforming the original regression problem (1) via a "multiscale" orthonormal linear transform $W$ called the Discrete Wavelet Transform (DWT), the following regression formulation is obtained:

$$Y_{j,k} = \mu_{j,k} + \varepsilon_{j,k}, \quad j = 0, \ldots, \log_2 n - 1, \quad k = 1, \ldots, 2^j,$$

and $k = 1$ for $j = -1$, where $j$ and $k$ are (respectively) scale and location parameters, $Y_{j,k}$ are the empirical wavelet coefficients of $X_t$, $\mu_{j,k}$ are the true wavelet coefficients of $\alpha(t/n)$ which need to be estimated, and $\varepsilon_{j,k}$ are the wavelet coefficients of $\varepsilon_t$. The sequence $\mu_{j,k}$ is often sparse, with most $\mu_{j,k}$'s being equal or close to zero, which motivates the use of simple thresholding techniques that do not estimate $\mu_{j,k}$ by zero if and only if the corresponding $Y_{j,k}$ exceeds a certain threshold in absolute value. This ensures that a large proportion of the noise $\varepsilon_{j,k}$ gets removed. The inverse DWT then yields an estimate of the original function $\alpha$.

The overwhelming majority of wavelet-based estimation techniques, such as those proposed by Donoho and Johnstone (1995), Johnstone and Silverman (2005a) or Abramovich et al. (2006), to name but a few, were devised under the assumption that the errors $(\varepsilon_t)_{t=1}^n$ formed an independent, identically distributed Gaussian sequence. This was partly due to the fact that in view of the orthonormality of $W$, the "wavelet noise" $\varepsilon_{j,k}$ was then also i.i.d. Gaussian, which facilitated both the choice of thresholds and the theoretical analysis of the resulting estimators. Johnstone and Silverman (1997) proposed an extension of the wavelet thresholding paradigm to stationary Gaussian noise.

In practice, the assumption of Gaussianity is violated in many important estimation problems. We list and discuss a selection of them below.

- *Poisson intensity estimation.* In Poisson intensity estimation, $X_t$ are modelled as independent $\text{Pois}\{\alpha(t/n)\}$ variables, which implies that $\varepsilon_t$ are centered Poisson. The mean and variance of $X_t$ are linked via the relationship $\text{var}(X_t) = h\{\mathbb{E}(X_t)\}$ with $h(u) = u$. This is in contrast to the i.i.d. Gaussian model in which $h(u) = \text{const}$. Recent examples of wavelet-based (or otherwise multiscale) Poisson intensity estimation techniques include the Bayesian methods of Kolaczyk (1999) and Timmermann and Nowak (1999), the multiscale likelihood technique of Kolaczyk and Nowak (2004), the Haar-Fisz method of Fryzlewicz and Nason (2004) and the extension of the latter proposed by Jansen (2006). Some of the above, as well as some other, techniques are reviewed in Besbeas et al. (2004). The work by Sardy et al. (2004), amongst other contributions, proposes an automatic smoothing procedure for Poisson data.

- *Nonparametric volatility estimation.* Nonparametric volatility estimation techniques are widely used in the finance industry (for example by *Risk-Metrics*, as detailed in their "Technical Document" available from http://www.riskmetrics.com/pdf/td4e.pdf). In this set-up, the $X_t$'s represent squared log-returns on a financial instrument and are modelled as independent and distributed as $X_t = \alpha(t/n) Z_t^2$, where $\mathbb{E}(Z_t^2) = 1$. Note that $\varepsilon_t = \alpha(t/n)(Z_t^2 - 1)$. Thus, the model is multiplicative and the variance function $h(u)$ is proportional to $u^2$. A multiscale Haar-Fisz technique for nonparametric volatility estimation was proposed by Fryzlewicz et al. (2006).

- *Spectral density estimation.* In spectral density estimation based on the periodogram, the $X_t$'s represent periodogram ordinates and are assumed to be asymptotically independent and asymptotically distributed as $\alpha(t/n) Z_t^2$, where $\alpha(t/n)$ represents the spectral density at frequency $t/n$, and $Z_t^2$ are Exp(1) random variables. This again renders the set-up multiplicative and, asymptotically, the variance function takes the form $h(u) = u^2$. Recent wavelet and multiscale approaches to periodogram smoothing include Moulin (1994), Neumann (1996), Gao (1997a), Pensky et al. (2007) and Fryzlewicz et al. (2008).

In the above examples, the variance function $h(u)$ is assumed to be known (as is the case in the work of Antoniadis and Sapatinas (2001) and Antoniadis et al. (2001)), and all of the multiscale approaches listed above, in one way or another, make use of its exact form in order to set the threshold at the "right" level: the threshold value usually depends on $\text{Var}(\varepsilon_{j,k})$, which involves $h$. However, in many estimation problems modelled by (1), it is clear that there exists a non-trivial mean-variance relatonship, but its exact form is unknown. Thus, it is not a priori clear what threshold values to use. For example, in gene expression data observed in microarray experiments, Rocke and Durbin (2001) identified that the variance of raw pixel intensities increased with their mean. An interesting mean-variance relationship arising in solar irradiance data was described in Fryzlewicz et al. (2007). Also, even if the variance function is "assumed to be known", the model which implies its particular form might have been chosen incorrectly. Thus, even in such a case it may often be safer to infer the form of the variance function from the data in order to set the thresholds at a suitable level.

A seemingly attractive solution to the problem of the unknown variance is to apply a data-driven variance stabilisation technique prior to smoothing the transformed data by means of a wavelet-based technique suitable for homoscedastic noise. Examples of data-driven variance-stabilising transforms include the ACE method of Breiman and Friedman (1985), the AVAS technique of Tibshirani (1988) as well as the method of Linton et al. (1997). In the context of one-colour microarray data, we mention the generalised log transformation proposed by Rocke and Durbin (2001) and the Spread-Versus-Level technique of Archer (2004). As these transforms operate on the original data, we refer to them as "time-domain" transforms below.

Despite its appealing modularity, a three-stage estimation procedure which involves a time-domain data-driven variance-stabilising transform followed by wavelet smoothing of the transformed data and finally the inverse transform, is not always recommendable. Firstly, the performance of time-domain data-driven variance stabilisation procedures is often less than satisfactory, as illustrated in Fryzlewicz et al. (2007). Secondly, due to the random and often highly nonlinear character of such transforms, theoretical properties of the resulting three-stage estimator may not be easy to establish.

Another possible route to follow is to treat the estimation problem (1) as an instance of the problem of estimating a function contaminated with locally stationary noise. Solutions to the latter were proposed, amongst others, by Gao (1997b) and von Sachs and MacGibbon (2000) (see also the references therein). However, adapted to our setting, these approaches would mean ignoring the fact that the local mean and variance were linked via a variance function $h$, and simply pre-estimating the evolution of the variance of $X_t$ over time. Thus, these methods can potentially be suboptimal in our context, as they do not take advantage of all available information.

In this paper, we propose an alternative approach to the problem of estimating $\alpha$ when $h$ exists but is unknown, which consists of first *estimating the variance function $h$*, and then *constructing a wavelet thresholding estimator of $\alpha$ which makes use of the estimate of $h$*. Our method, termed the *data-driven wavelet-Fisz* estimation technique, overcomes the drawbacks mentioned above in that (a) it performs well in practice, (b) it only requires that the variance function $h$, and not the target function $\alpha$, be "smooth", and (c) the theoretical performance of the final estimator of $\alpha$ is possible to quantify and near-optimal. In addition, our estimator of $\alpha$ is rapidly computable and easy to implement. Its simple modification can also be used in cases in which the variance function $h$ is known. Thus, it is applicable to a wide range of problems, including all of the examples mentioned above. The added benefit of our approach is that as well as estimating $\alpha$, it also returns an estimate of $h$, which may be of interest to the analyst.

Our estimator of $h$ is a simple Nadaraya-Watson estimator, and is inspired by (but simpler than) the variance function estimator used by Chiou and Müller (1999) in a different context. The reason for preferring simplicity here is that in order to derive theoretical properties of the resulting estimator of $\alpha$, we need to establish an exponential inequality for the estimator of $h$. The latter piece of theory forms a large part of this work, and we hope that it may be of independent interest. We note that a large deviation theory for a class of Nadaraya-Watson estimators was obtained by Louani (1999) and Joutard (2006). However, it was done in a simple nonparametric regression set-up with independent errors, which is not applicable in the context of variance function estimation. We also note that the theoretical part of our work does not address the automatic choice of the smoothing parameter occurring in the estimator of $h$. However, the simulations section provides detailed practical recommendations as to the choice of this parameter. This selection is also performed automatically in the software package which accompanies this paper. Finally, we note that what we mean by the term "variance function estimation" is different from the use of the same term by some other authors, for example Cai and Wang (2007) and Wang et al. (2008), who consider the estimation of the variance function as a function of *time*, not as a function of the *mean*. However, a similar element in our work and the above-cited articles is the fact that in both estimation problems, a preliminary estimator of the mean is used to estimate the variance function.

It is interesting to note that the algorithm for computing our estimator of $\alpha$ can be decomposed into three separate stages, the first and last of which is a particular data-driven variance-stabilising transform, and its inverse, respectively. In the rare cases when the exact standard deviation of the noise is known, the simplest *time-domain* variance stabilisation procedure consists of using this quantity to pre-divide the original regression set-up coordinate-wise, as discussed by Gao (1997b). Unlike this and other existing transforms, some of which are listed above, our variance-stabilising transform is *performed in the wavelet domain*, as opposed to the time domain. Roughly speaking, the transform consists of dividing each empirical wavelet coefficient $Y_{j,k}$ by an estimate of its own standard deviation, the latter involving the estimate of $h$. In a non-wavelet context, similar ratio statistics (for a known function $h$) were studied by Fisz (1955), which justifies the name of our procedure.

We also mention that this paper was inspired by our earlier work Motakis et al. (2006) and Fryzlewicz et al. (2007), where we proposed computational procedures related to that described here. However, they were not accompanied by any theoretical analysis, partly because any such analysis appeared challenging due to the level of complexity of the proposed algorithms. Indeed, one of the aims of this paper is to provide a procedure which is simple enough to be theoretically tractable, but also performs well in practice.

The paper is organised as follows. In Section 2, we describe our model and estimation problem. In Section 3, we introduce and analyse our wavelet-Fisz technique in the case when the variance function $h$ is known. Section 4 describes the Nadaraya-Watson estimator of $h$ and considers its theoretical properties. In Section 5, we show how the estimator of $h$ from Section 4 is used in our *data-driven wavelet-Fisz estimator* of $\alpha$, and establish the mean-square near-optimality of the latter. In Sections 6 and 7, we discuss various implementation aspects of our estimators of $h$ and $\alpha$, respectively. Section 8 demonstrates how our wavelet estimator leads to a new data-driven variance-stabilising transform performed in the wavelet-domain. Section 9 concludes, and the proofs of our results are in three appendices. R code implementing our method has been made available on http://www.maths.bris.ac.uk/~mapzf/ddwf/ddwf.html

## 2. Model and preliminaries

We consider the regression model

$$X_t = \alpha(t/n) + \varepsilon_t, \quad t = 1, \ldots, n, \tag{2}$$

where $X_1, \ldots, X_n$ are assumed to be nonnegative and independent, with $\mathbb{E}(X_t) = \alpha(t/n)$ and $\mathrm{Var}(X_t) = \mathrm{Var}(\varepsilon_t) = h\{\alpha(t/n)\}$. Our task is to estimate $\alpha$ nonparametrically via nonlinear wavelet shrinkage, assuming that the function $h$ is not necessarily known.

We place the following assumption on the function $\alpha$.

**Assumption 1.** *For the function $\alpha(z) : [0,1] \mapsto \mathbb{R}$, we denote $\underline{\alpha} = \inf_z \alpha(z)$ and $\overline{\alpha} = \sup_z \alpha(z)$. We assume*

    *(i) $\alpha$ is of finite total variation over $[0,1]$,*
    *(ii) $0 < \underline{\alpha} \leq \overline{\alpha} < \infty$.*

Assumption 1(i) is a mild smoothness assumption on $\alpha$. Since the class of bounded variation functions also includes functions of a low degree of regularity, the choice of nonlinear wavelet shrinkage as the preferred estimation method appears natural in this context.

As mentioned earlier, our estimator of $\alpha$ can be viewed as using the principle of variance stabilisation (in the wavelet domain). Many time-domain variance-stabilising transforms, such as the square-root transform for Poisson data or the log transform for scaled $\chi^2$ data, would require that the function $\alpha$ be bounded from below, as specified in Assumption 1(ii). Therefore, it is not surprising that we also require this assumption to hold.

Further, we impose the following assumption on $h$.

**Assumption 2.** *For the function $h : [0,\infty) \mapsto [0,\infty)$ we denote $\underline{h} = \inf_{u \in [\underline{\alpha}, \overline{\alpha}]} h(u)$ and $\overline{h} = \sup_{u \in [\underline{\alpha}, \overline{\alpha}]} h(u)$. We assume*

    *(i) $0 < \underline{h} \leq \overline{h} < \infty$,*
    *(ii) $h$ is non-decreasing,*
    *(iii) $h$ is Lipschitz-continuous of order 1 on $u \in [\underline{\alpha}, \overline{\alpha}]$ with constant $H$,*
    *(iv) there exist $\tilde{\delta}, \bar{\delta}, \tilde{H} > 0$ such that $h^{\bar{\delta}}$ is Hölder-continuous on $u \in [0,\infty)$ with Hölder exponent $\bar{\delta}$ and constant $\tilde{H}$.*

The class of distributions with a variance function $h$ satisfying Assumption 2 includes, amongst others, the Poisson distribution, for which $h(u) = u$, and distributions of the form $X_t = \alpha(t/n)Z_t$ where $\{Z_t\}_t$ are i.i.d. with $\mathbb{E}(Z_t) = 1$, for which $h(u)$ is proportional to $u^2$.

Finally, we make the following assumption about the central moments of $\varepsilon_t$.

**Assumption 3.** *We assume that there exists a positive constant $K$ and a non-negative constant $\gamma$ such that*

$$\mathbb{E}|\varepsilon_t|^l = \mathbb{E}|X_t - \alpha(t/n)|^l \leq (l!)^{1+\gamma} K^{l-2} h\{\alpha(t/n)\}$$

*for $l = 3, 4, \ldots$ and all $t$.*

Assumption 3 is natural and common in the context of wavelet estimation in non-Gaussian noise, see for example Neumann (1996). It is satisfied by many standard distributions, including, amongst others, Poisson and gamma. Roughly speaking, it ensures that local sums of $X_t$ are asymptotically normal in a certain asymptotic regime and in a certain required "strong" sense.

## 3. Wavelet-Fisz estimation for *h* known

In this section, we aim to estimate $\alpha$ using nonlinear wavelet shrinkage assuming that the variance function $h$ is known. Throughout the paper, we assume basic familiarity with the Discrete Wavelet Transform (DWT). We refer the reader to Mallat (1989) for a description of the DWT, and to Vidakovic (1999) for an excellent overview of wavelet methods in statistics.

We now describe, step by step, our algorithm for computing the wavelet-Fisz estimator of $\alpha$ if the function $h$ is known.

1. The starting point is the DWT of the observed data $\{X_t\}_{t=1}^n$ with respect to an orthonormal basis of compactly supported wavelets. Later, in Assumption 4, we will specify additional technical conditions on the wavelets. The DWT converts the regression problem (1) into a regression problem in the wavelet domain

$$Y_{j,k} = \mu_{j,k} + \varepsilon_{j,k}, \quad j = 0, \ldots, J-1, \quad k = 1, \ldots, 2^j,$$

where $J = \log_2 n$, with the only "smooth" coefficient indexed by $(j,k) = (-1, 1)$. The variables $Y_{j,k}$ are the empirical wavelet coefficients of $X_t$, the constants $\mu_{j,k}$ are the wavelet coefficients of $\alpha(t/n)$ which need to be estimated, and the "wavelet noise" variables $\varepsilon_{j,k}$ are the wavelet coefficients of $\varepsilon_t$. The sequence $\mu_{j,k}$ will often be sparse, with most $\mu_{j,k}$'s being equal or close to zero.

2. We then separate the indices $(j,k)$ into two groups: those corresponding to the coarser scales $0 \le j \le J^* - 1$, for which $\varepsilon_{j,k}$ will be asymptotically normal, and those corresponding to the finer scales $J^* \le j \le J - 1$, which will be "ignored" in the estimation procedure. To be more precise, we define $J^*$ and a set $\mathcal{I}_n$ as follows:

$$\mathcal{I}_n = \{(j,k) \in \mathbb{Z}^2 \quad : \quad 1 \le k \le 2^j; \quad 0 \le j \le J^* - 1; \quad 2^{J^*} = n^{1-\epsilon}\},$$

for a fixed $\epsilon \in (0, 1/3]$. The choice of $1/3$ as the upper bound for $\epsilon$ is linked to the postulated smoothness of $\alpha$. The reader is referred to Section 3.3 of Fryzlewicz et al. (2008) for a more detailed discussion of this issue.

3. As the sequence $\mu_{j,k}$ is likely to be sparse, with most $\mu_{j,k}$'s being equal or close to zero, we use a simple thresholding technique, which does not estimate $\mu_{j,k}$ by zero if and only if the corresponding $Y_{j,k}$ exceeds a certain threshold in absolute value. This ensures that a large proportion of the noise $\varepsilon_{j,k}$ gets removed.

In wavelet function estimation with Gaussian errors, possibly the simplest ("universal") threshold, advocated by Donoho and Johnstone (1994), takes the form

$$\lambda_{j,k}^U = \{2 \operatorname{Var}(\varepsilon_{j,k}) \log(\#\mathcal{I}_n)\}^{1/2}, \tag{3}$$

where $\#A$ is the cardinality of the set $A$. Since in the set $\mathcal{I}_n$, our wavelet coefficients $\varepsilon_{j,k}$ are *asymptotically* Gaussian, we wish to explore the possibility of applying an analogous threshold in our set-up. To effect this idea, we need to determine $\mathrm{Var}(\varepsilon_{j,k})$. Denoting by $\{\psi_{j,\tau}\}_\tau$ the elements of the discrete wavelet vector at scale $j$, we find

$$\mathrm{Var}(\varepsilon_{j,k}) = \mathrm{Var}\left(\sum_t \psi_{j,k-t}\varepsilon_t\right) = \sum_t \psi_{j,k-t}^2 h\{\alpha(t/n)\} \qquad (4)$$

Obviously, $\alpha(t/n)$ is unknown and needs to be pre-estimated. For simplicity and speed of computation, we use the same pre-estimate for each $t$ in $\mathrm{supp}(\psi_{j,k-\cdot})$, namely the following localised mean of $X_t$:

$$\widehat{\alpha(t/n)} = \sum_q \kappa_{j,k-q}X_q, \qquad (5)$$

where the constants $\kappa_{j,\tau}$ satisfy Assumption 5 below. As the discrete wavelet vectors are normalised so that $\sum_t \psi_{j,k-t}^2 = 1$, we obtain our "estimated" thresholds as

$$\hat{\lambda}_{j,k} = h^{1/2}\left(\sum_q \kappa_{j,k-q}X_q\right)\{2\log(\#\mathcal{I}_n)\}^{1/2}. \qquad (6)$$

(As a side remark, we note that in the absence of an assumed mean-variance relationship, Gao (1997b) estimates thresholds as proportional to the "running median absolute deviation" estimate of the standard deviation of the heteroscedastic noise.)

We use our estimated thresholds to estimate $\mu_{j,k}$ in the set $\mathcal{I}_n$ via either the soft or the hard thresholding rule:

$$\hat{\mu}_{j,k}^{(s)} = \mathrm{sign}(Y_{j,k})\max\left(|Y_{j,k}| - \hat{\lambda}_{j,k}, 0\right) \qquad (7)$$

$$\hat{\mu}_{j,k}^{(h)} = Y_{j,k}\mathbb{I}\left(|Y_{j,k}| \geq \hat{\lambda}_{j,k}\right), \qquad (8)$$

where $\mathbb{I}(\cdot)$ is the indicator function. Outside the set $\mathcal{I}_n$, we simply estimate $\mu_{j,k}$ by zero, that is

$$\hat{\mu}_{j,k}^{(s)} = \hat{\mu}_{j,k}^{(h)} = 0 \quad \text{for} \quad (j,k) \in \{j \geq 0\} \cap \mathcal{I}_n^c.$$

4. Let $(e)$ denote either one of: $(s)$ or $(h)$. The inverse DWT of the sequence $\hat{\mu}_{j,k}^{(e)}$ yields our wavelet-Fisz estimator $\tilde{\alpha}^{(e)}$.

A few remarks are in order.

*Stability.* Let $\mathrm{stdev}(X)$ denote the standard deviation of a random variable $X$. Looking back at the derivation in formula (4), another "obvious" way of estimating $\mathrm{stdev}(\varepsilon_{j,k})$ would be to set $\widetilde{\mathrm{stdev}(\varepsilon_{j,k})} = \left\{\sum_t \psi_{j,k-t}^2 h\{X_t\}\right\}^{1/2}$. However, comparing it to our estimator $\widehat{\mathrm{stdev}(\varepsilon_{j,k})} = h^{1/2}\left(\sum_t \kappa_{j,k-t}X_t\right)$ from formula (6), it is easily seen that the latter typically involves lower powers of $X_t$, and thus is potentially a more "stable" statistic. As an example, consider $h(u) = u^2$.

In this case, $\widetilde{\text{stdev}(\varepsilon_{j,k})}$ is a localised $l_2$ norm of $X_t$, whereas $\widehat{\text{stdev}(\varepsilon_{j,k})}$ is a localised $l_1$ norm. A similar comment applies in the case $h(u) \sim u^\beta$ for all $\beta > 1$.

*Link to maximum likelihood estimation.* If $\alpha(t/n)$ is constant over the support of $\kappa_{j,k-\cdot}$, then, by the invariance principle of maximum likelihood estimators, our estimator $\widehat{\text{stdev}(\varepsilon_{j,k})}$ is precisely the maximum likelihood estimator of $\text{stdev}(\varepsilon_{j,k})$, provided that $\kappa_{j,k-t} = \text{const for } t \in \text{supp}(\kappa_{j,k-\cdot})$.

*The name "wavelet-Fisz".* Note that the argument of the indicator function in (8) can be rewritten as

$$\left| \frac{\sum_t \psi_{j,k-t} X_t}{h^{1/2} \left( \sum_t \kappa_{j,k-t} X_t \right)} \right| \geq \{2 \log(\#\mathcal{I}_n)\}^{1/2} \tag{9}$$

In a non-wavelet context, ratio transformations similar to that on the left-hand side of (9) and the asymptotic normality of the resulting random variables were studied by Fisz (1955), which justifies the name of our procedure. The division in (9) provides a degree of "variance stabilisation": note that the threshold $\{2 \log(\#\mathcal{I}_n)\}^{1/2}$ is suitable for standard homoscedastic normal noise. In this sense, our procedure can be viewed as being based on the principle of variance stabilisation in the wavelet domain. We expand on this issue later in Section 8.

*Link to Fryzlewicz et al. (2008).* We note that in the special case $h(u) = u^2$, our estimation algorithm is equivalent to the method proposed by Fryzlewicz et al. (2008) in the context of spectral density estimation.

We now establish the mean-square convergence rate of our wavelet-Fisz estimator $\tilde{\alpha}^{(e)}$. In order to do so, we specify assumptions on the wavelets $\psi_{j,k}$ and the constants $\kappa_{j,k}$.

**Assumption 4.** *The discrete wavelets used in the construction of $\tilde{\alpha}^{(e)}$ are derived from a continuous-time orthonormal wavelet basis of $L_2[0,1]$, $\{\phi_{0,k}(z)\}_k \cup \{\psi_{j,k}(z)\}_{j \geq 0,k}$, where $\phi_{j,k}(z) = 2^{j/2}\phi(2^j z - k)$ and $\psi_{j,k}(z) = 2^{j/2}\psi(2^j z - k)$. The "mother" and "father" wavelet functions $\psi$ and $\phi$ are assumed to satisfy, for some $r > m$ (with $m$ given in Theorem 1 below),*

*(i) $\phi$ and $\psi$ are in the space $C^r$,*
*(ii) $\int \phi(z)dz = 1$,*
*(iii) $\int \psi(z)z^l dz = 0$ for all $0 \leq l \leq r$.*

Assumption 4 defines the so-called $r$-regularity of the wavelet basis, and is commonly used in wavelet function estimation.

**Assumption 5.** *The constants $\kappa_{j,\tau} \geq 0$ are such that*

$$\sum_\tau \kappa_{j,\tau} = 1$$

$$\sum_\tau \kappa_{j,\tau}^2 = O(2^{j-J})$$

$$\max_{\tau} \kappa_{j,\tau} = O(2^{j-J})$$
$$\operatorname{supp} \kappa_{j,\cdot} = \operatorname{supp} \psi_{j,\cdot}$$

*for all* $0 \le j \le J^* - 1$.

Note that each of the vectors $\kappa_j$ can be interpreted as a linear filter which computes the local mean of $X_t$ over the support of the vector $\psi_j$ in (5).

As has now become standard in the wavelet literature, we assume that the unknown function $\alpha$ is in a Besov ball of radius $C > 0$ on $[0,1]$, $\mathcal{F}^m = \mathcal{F}^m_{p,q}(C)$, where $m > 0$ and $0 < p, q \le \infty$. Roughly speaking, the not necessarily integer parameter $m$ indicates the number of derivatives of $\alpha$, where their existence in required in the $L_p$-sense, and thus $p$ can be viewed as a measure of the inhomogeneity of $\alpha$. The additional parameter $q$ provides a further finer gradation. Besov classes include the traditional Hölder and Sobolev classes ($p = q = \infty$ and $p = q = 2$, respectively). If the regularity $r$ of a wavelet basis is greater than $m$ and if other conditions of Assumption 4 hold, then the membership of $\alpha$ in $\mathcal{F}^m$ can be characterised in terms of the wavelet coefficients $\mu'_{j,k} = \mu_{j,k} n^{-1/2}$ of the function $\alpha$ in the following way. Define the Besov sequence ball of radius $C$ as

$$b^m_{p,q}(C) = \left\{ \mu'_{j,k} : \sum_{j \ge 0} 2^{jsq} \|\mu'_j\|^q_p \le C^q \right\},$$

where $s = m + 1/2 - 1/p$ and $\|\mu'_j\|^p_p = \sum_{k=1}^{2^j} |\mu'_{j,k}|^p$. If Assumption 4 holds, then $\alpha$ is in $\mathcal{F}^m$ if and only if $\{\mu'_{j,k}\}_{j,k}$ is in $b^m_{p,q}(C)$. The reader is referred to Meyer (1992) for rigorous definitions and a detailed study of Besov spaces.

Denote $\|v\|^2_{L_2[0,1]} = \int_0^1 |v(u)|^2 du$. We are now ready to state a result on the mean-square rate of convergence of our wavelet-Fisz estimator $\tilde{\alpha}^{(e)}$.

**Theorem 1.** *Let* (e) *denote either one of:* (s) *or* (h). *Suppose that Assumptions 1, 2, 3, 4 and 5 hold. We have*

$$\sup_{\alpha \in \mathcal{F}^m} \mathbb{E}\|\tilde{\alpha}^{(e)} - \alpha\|^2_{L_2[0,1]} = O\left\{ (\log n/n)^{2m/(2m+1)} \right\}.$$

The rate $O\left\{ n^{-2m/(2m+1)} \right\}$ is the best possible mean-square error rate for Besov spaces, and our wavelet-Fisz estimator achieves it up to a logarithmic term, attaining the same rate as the universal thresholding estimator in the case of i.i.d. Gaussian noise. We mention that linear estimators, such as kernel estimators, cannot attain the optimal mean-square error rate (not even up to a logarithmic factor) for $p < 2$.

## 4. Estimation of the variance function $h$

In this section, we assume that the function $h$ is unknown, and we propose to estimate it by means of a Nadaraya-Watson estimator $\hat{h}$. Later, in Section 5, $\hat{h}$ will be used in the data-driven wavelet-Fisz estimator of $\alpha$. In order to establish the mean-square convergence of the latter estimator, we need to be able to determine large deviation properties of $\hat{h}$. Indeed, the main aim of this section is to demonstrate an exponential inequality for $\hat{h}$ (which will be sufficient for our purposes).

The estimator $\hat{h}$ is constructed as follows. We start with a preliminary estimator of $\alpha(t/n)$ defined by

$$\hat{\alpha}_t = \frac{1}{2M+1} \sum_{p=t-M}^{t+M} X_p, \qquad (10)$$

where the choice of $M$ will be discussed in the paragraph underneath Theorem 2, as well as in Section 6. We define empirical residuals from this fit by $\hat{\varepsilon}_t = X_t - \hat{\alpha}_t$. Our Nadaraya-Watson estimator $\hat{h}$ performs kernel smoothing of the squared empirical residuals $\hat{\varepsilon}_t^2$. More specifically, we use a kernel function $K$ which satisfies the following assumption.

**Assumption 6.** *The function $K : [-1/2, 1/2] \mapsto \mathbb{R}$ is nonnegative, bounded, integrates to one and is Lipschitz-continuous of order 1 with constant $L$. We denote $\check{K} = \max_z K(z)$.*

We define

$$\begin{aligned} W_{nt}(u) &= \frac{1}{nb} K\left(\frac{\alpha(t/n) - u}{b}\right) \\ \hat{W}_{nt}(u) &= \frac{1}{nb} K\left(\frac{\hat{\alpha}_t - u}{b}\right), \end{aligned} \qquad (11)$$

where the choice of $b$ will also be discussed in the paragraph underneath Theorem 2, as well as in Section 6. The Nadaraya-Watson estimator $\hat{h}$ is defined by

$$\hat{h}(u) = \frac{\sum_{t=1}^{n} \hat{W}_{nt}(u) \hat{\varepsilon}_t^2}{\sum_{t=1}^{n} \hat{W}_{nt}(u)}.$$

We now list and clarify a number of assumptions which will be used in proving the main result of this section.

**Assumption 7.** *We have $\mathrm{Var}(\varepsilon_t^2) \geq \underline{c} > 0$, uniformly over $t$.*

**Assumption 8.** *Denote $Z_t = |\hat{\alpha}_t - \mathbb{E}(\hat{\alpha}_t)|$. We assume that there exists a positive constant $C_2$ such that*

$$\mathrm{Var}(\hat{\alpha}_t) = \mathrm{Var}(\hat{\alpha}_t - \mathbb{E}(\hat{\alpha}_t)) \leq C_2 \, \mathrm{Var}(Z_t),$$

*uniformly over $t$.*

**Assumption 9.** *We assume that there exist positive constants $a$, $d$ such that*

$$P\left(\sup_{t=1,\ldots,n} |\varepsilon_t| \leq a \, \log^d n\right) = 1 - O(n^{-2}).$$

**Assumption 10.** *There exists a function $c(u)$ such that*

$$0 < c_1 \leq c(u) \leq \sum_{t=1}^{n} W_{nt}(u)$$

*uniformly over $n$ and $b$, for all $u \in \text{range}\{\alpha(z)\}$.*

**Assumption 11.** *Let the constants $\kappa_{j,\tau}$ satisfy Assumption 5. We assume that the function $\alpha(z)$ is such that for all $(j,k) \in \mathcal{I}_n$,*

$$\sum_q \kappa_{j,k-q}\alpha(q/n) \in \text{range}\{\alpha(z)\}.$$

Assumption 7 ensures that $\varepsilon_t^2$ is a non-degenerate random variable.

For any random variable $Y$, it is easy to see that $\text{Var}|Y| \leq \text{Var}(Y)$. Assumption 8 guarantees that the converse is true, up to a constant. This is not a restrictive assumption: we expect that $\hat{\alpha}_t - \mathbb{E}(\hat{\alpha}_t)$ will be close to $N(0, \sigma^2)$, and for $Y \sim N(0, \sigma^2)$ we have $\text{Var}(Y) = \frac{\pi}{\pi-2}\text{Var}|Y|$.

Assumption 9 is satisfied for all distributions whose tail decays exponentially.

We now demonstrate that Assumption 10 is satisfied for functions $\alpha(z)$ which are piecewise Lipschitz-continuous of order 1 with a finite number of breakpoints. For clarity, we only show it for the triangular kernel $K(z) = (-4|z| + 2)\mathbb{I}(|z| \leq 1/2)$, but the proof remains almost unchanged for other kernels.

**Proposition 1.** *Let $K(v) = (-4|v| + 2)\mathbb{I}(|v| \leq 1/2)$ and let $\alpha(z)$ be piecewise Lipschitz-continuous of order 1 with a finite number of jumps. There exists a positive constant $c_1$ such that $\sum_{t=1}^{n} W_{nt}(u) \geq c_1$ uniformly over $n$, $b$, and $u \in \text{range}\{\alpha(z)\}$.*

We note that Assumption 10 can be relaxed to include functions which are piecewise Hölder continuous, at the expense of worse rates of tail decay in Theorem 2. As an interesting example, we note that Donoho and Johnstone's (1994) benchmark signals *blocks*, *bumps*, *doppler* and *heavisine* are all piecewise Lipschitz-continuous of order 1.

Assumption 11 ensures that "local averages" of the function $\alpha$ lie within the range of $\alpha$.

We now state an exponential inequality for the estimator $\hat{h}(u)$, which is the main result of this section.

**Theorem 2.** *Suppose that Assumptions 1–3, 6–11 hold, and that the constants $\kappa_{j,\tau}$ satisfy Assumption 5. Let $b_n$, $d_n$ and $i_n$ be any fixed sequences such that $b_n = o\left\{(n/M)^{1/(6+4\gamma)}\right\}$, $d_n = o\left\{\min_j 2^{\frac{J-j}{2(1+\max(\gamma,1))}}\right\}$ and $i_n = o\left\{(nb^2)^{1/(10+12\gamma)}\right\}$, where $M$, $b$ and $\gamma$ are defined in formulae (10) and (11) and Assumption 3, respectively; $2^J = n$, and the range of $j$ is $0 \leq j \leq J^* - 1$. Let $\delta_1$ be any positive quantity such that $\delta_1 < c_1$ where $c_1$ is defined in Assumption 10, and define $\delta' = \delta(c_1 - \delta_1)/2$, where $\delta$ appears in the exponential inequality below. Assume*

$$\delta_1 n^{1/2}b^2 - Mn^{-1/2} - n^{1/2}M^{-1/2} \;\; \rightarrow \;\; \infty \qquad (12)$$

$$\delta'n^{1/2}b^2 \log^{-2d} n - Mn^{-1/2} - n^{1/2}M^{-1/2} \;\; \rightarrow \;\; \infty \qquad (13)$$

$$\delta \;\; > \;\; bH, \qquad (14)$$

*as $M, n/M, nb \rightarrow \infty$ and $b \rightarrow 0$, where $d$ is defined in Assumption 9 and $H$ is the Lipschitz constant for $h(u)$ over $u \in \text{range}\{\alpha(z)\}$. In the asymptotic limit, as $M, n/M, nb \rightarrow \infty$ and $b \rightarrow 0$, we have, uniformly over $(j, k) \in \mathcal{I}_n$,*

$$P\left(\left|\hat{h}\left(\sum_q \kappa_{j,k-q} X_q\right) - h\left(\sum_q \kappa_{j,k-q}\alpha(q/n)\right)\right| \geq \delta\right)$$
$$\leq C_3(2M+1)\left\{4 - \Phi(\min(a_n, b_n)) - 2\Phi(\min(e_n, b_n)) - \Phi(\min(g_n, b_n))\right\}$$
$$+C_4\left\{3 - \Phi(\min(c_n, d_n)) - 2\Phi(\min(f_n, d_n))\right\} + C_5\left\{1 - \Phi(\min(h_n, i_n))\right\}$$
$$+O(n^{-2}),$$

*where $C_3$, $C_4$ and $C_5$ are positive constants, $\Phi$ is the cdf of the standard normal, and*

$$
\begin{aligned}
a_n &= O(\delta_1 n^{1/2}b^2 - Mn^{-1/2} - n^{1/2}M^{-1/2})\\
c_n &= O(\delta_1 b^2 \min_j 2^{(J-j)/2})\\
e_n &= O(\delta'n^{1/2}b^2 \log^{-2d} n - Mn^{-1/2} - n^{1/2}M^{-1/2})\\
f_n &= O(\delta'b^2 \min_j 2^{(J-j)/2} \log^{-2d} n)\\
g_n &= O(\delta'n^{1/2}b \log^{-d} n - Mn^{-1/2} - n^{1/2}M^{-1/2})\\
h_n &= O\{(\delta - bH)n^{1/2}b\}.
\end{aligned}
$$

*Explanation of the rates.* We take $M = O(n^\vartheta)$ for $\vartheta \in (0, 1)$ and $b = O(n^{-\zeta})$ for $\zeta \in (0, 1)$ and investigate conditions on $\vartheta$ and $\zeta$ which ensure that the assumptions of Theorem 2 are satisfied. Clearly, if $\zeta < 1/2$, then $b_n$, $d_n$ and $i_n$ can be chosen such that they are all of order $O(n^\varsigma)$, for $\varsigma > 0$. Fixing $\delta_1$ to be constant, and assuming that $\delta$ and $\delta'$ either are constant or tend to zero no faster than logarithmically in $n$, we have that conditions (12) – (14) are satisfied if both $1 > 2\zeta + \vartheta$ and $\vartheta/4 > \zeta$. Finally, to ensure that the sequences $a_n$, $c_n$, $e_n$, $f_n$, $g_n$ and $h_n$ are all of order $O(n^\varsigma)$ for $\varsigma > 0$, we additionally require that $\epsilon/4 > \zeta$, where $2^{J^*} = n^{1-\epsilon}$. Thus, in the $(\epsilon, \vartheta, \zeta)$ space, the set $A$ of parameter configurations which are "admissible" in the above sense has the form

$$A = \{(\epsilon, \vartheta, \zeta) : \epsilon \in (0, 1/3], \zeta \in (0, \min(\vartheta/4, 1/2 - \vartheta/2, \epsilon/4))\}. \qquad (15)$$

In view of the above discussion, the following corollary can be formulated.

**Corollary 1.** *Suppose that Assumptions 1, 2, 3, 5, 6, 7, 8, 9, 10 and 11 hold. Let $A$ be as defined in (15), and let $M = O(n^\vartheta)$ and $b = O(n^{-\zeta})$. If $(\epsilon, \vartheta, \zeta) \in A$ and $\delta \geq O(\log^v n)$ for some $v < 0$, then there exists $\varsigma > 0$ such that, uniformly* *over* $(j,k) \in \mathcal{I}_n$,

$$P\left( \left| \hat{h}\left( \sum_q \kappa_{j,k-q} X_q \right) - h\left( \sum_q \kappa_{j,k-q} \alpha(q/n) \right) \right| \geq \delta \right)$$
$$\leq O\left\{ n^{\vartheta} \exp\left( -\frac{n^{2\varsigma}}{2} \right) + n^{-2} \right\}.$$

## 5. Data-driven wavelet-Fisz estimation (for $h$ unknown)

Our algorithm for computing the wavelet-Fisz estimator of $\alpha$ if the function $h$ is unknown (which we also call the *data-driven* wavelet-Fisz estimator of $\alpha$), proceeds in the same way as the wavelet-Fisz algorithm for $h$ known, described in detail in Section 3, the only exception being that we use our estimate $\hat{h}$ from Section 4, instead of the true $h$.

To be more precise, we replace our thresholds $\hat{\lambda}_{j,k}$, defined in formula (6) and subsequently used in the threshold estimators $\hat{\mu}_{j,k}^{(s)}$ and $\hat{\mu}_{j,k}^{(h)}$ (see formulae (7) and (8)), with thresholds

$$\tilde{\lambda}_{j,k} = \hat{h}^{1/2}\left( \sum_q \kappa_{j,k-q} X_q \right) \{2 \log(\#\mathcal{I}_n)\}^{1/2}. \tag{16}$$

We denote the thus constructed data-driven wavelet-Fisz estimator of $\alpha$ by $\bar{\alpha}^{(e)}$, where $(e)$ is one of: $(s)$ (soft thresholding) and $(h)$ (hard thresholding). The following theorem quanitifies the mean-square rate of convergence of $\bar{\alpha}^{(e)}$.

**Theorem 3.** *Let $(e)$ denote either one of: $(s)$ or $(h)$. Suppose that Assumptions 1–11 hold. Let $A$ be as defined in (15), and let $M = O(n^{\vartheta})$ and $b = O(n^{-\zeta})$. If $(\epsilon, \vartheta, \zeta) \in A$, then*

$$\sup_{\alpha \in \mathcal{F}^m} \mathbb{E}\|\bar{\alpha}^{(e)} - \alpha\|_{L_2[0,1]}^2 = O\left\{ (\log n/n)^{2m/(2m+1)} \right\}.$$

Comparing this result with Theorem 1, we note that the estimator $\bar{\alpha}^{(e)}$ does not exhibit any loss of asymptotic efficiency compared to $\tilde{\alpha}^{(e)}$ despite the fact that it uses an estimate of $h$ rather than the true $h$.

## 6. Estimation of $h$ − implementation issues

This section briefly describes the outcome of an extensive simulation study aimed at assessing the performance of the estimator $\hat{h}(u)$ for various parameter configurations. Recall that $h(u)$ is assumed to be non-decreasing (see Assumption 2). As $\hat{h}(u)$ is not guaranteed to be non-decreasing, in practice we used the following computational "correction" to $\hat{h}(u)$. Having obtained $\hat{h}(u)$, we input it into the (automatic) "pool-adjacent-violators" algorithm for least-squares isotone regression, described in detail in Johnstone and Silverman (2005b), Section 6.3. The resulting estimate, denoted hereafter by $\dot{h}(u)$, is a non-decreasing, piecewise constant function of $u$, which is as close as possible to $\hat{h}(u)$ in the least-squares sense. Empirically, we found that the use of $\dot{h}(u)$, rather than $\hat{h}(u)$, in the estimator $\bar{\alpha}^{(e)}$, results in a slightly superior performance of the latter. (We note that literature on monotone estimation of variance function is sparse; Dette and Pilz (2007) consider estimation of a monotone conditional variance in nonparametric regression.)

Having investigated various choices of the span $M$ and the bandwidth $b$ for a range of test functions and noise distributions, we found that $\dot{h}(u)$ performed particularly well for "small" values of $M$. Any value of $M \leq 3$ consistently resulted in good estimates. The examples later in this section use the value $M = 3$. The estimator seems to be less sensitive to the choice of $b$: this is due to the computational correction (described above), which "smooths out" any remaining wiggles in $\hat{h}(u)$. We recommend an "automatic" choice of $b$, such as that described in Gasser et al. (1991) and conveniently implemented in the routine `glkerns` from the R package `lokern`. We also use the default kernel function $K(\cdot)$ from the above routine.

We briefly illustrate the performance of $\dot{h}(u)$ on 4 simulated datasets. The models are: the Poisson model, whereby $X_t \sim \text{Pois}\{\alpha(t/n)\}$, and the exponential model, in which $X_t \sim \alpha(t/n)\,\text{Exp}(1)$. With each model, we use two functions $\alpha(z)$: the *blocks* function, scaled and shifted to have the minimum (maximum) value of 1 (22.6), and the *bumps* function, with the minimum (maximum) value of 3 (23.21). Both functions are sampled at 2048 equispaced points.

Figure 1 shows sample paths from each model, together with the corresponding estimates $\dot{h}^{1/2}(u)$. The estimator performs well in all cases. The good performance is not incidental: indeed, we found that for the parameter choices described above, the estimator $\dot{h}(u)$ performed well across all simulated examples.

## 7. Estimation of $\alpha$ – implementation issues

This section discusses the choice of the various parameters for our data-driven wavelet-Fisz estimator $\bar{\alpha}^{(e)}$. The examples in this section use Haar wavelets: this choice is motivated in Section 8 below. As a default option, we use translation-invariant (see Nason and Silverman (1995)) hard thresholding with $J^* = J - 2$, as this parameter configuration seems to offer the best empirical performance. We use the variance estimator $\dot{h}(u)$ described in Section 6 with $M = 1$. For each $j, k$, we choose the parameters $\kappa_{j,k-t}$ to be constant for $t \in \text{supp}(\psi_{j,k-\cdot})$. This is a natural choice for Haar wavelets as the coefficients $\sum_q \kappa_{j,k-q} X_q$ are already available to us as "by-products" of the discrete Haar transform.

Figure 2 shows the outcome of our estimation procedure described above for the sample paths from Figure 1. It is clear that our procedure performs very well for Poisson noise. Performance for exponential noise is also satisfactory given how noisy the original signals are: indeed, it is extremely hard to identify some of the features of the clean *bumps* and *blocks* signals from the visual inspection of the corresponding exponential datasets. We mention again that our estimation procedure "does not know" any characteristics of the noise, and estimates the variance function $h(u)$ from the data.

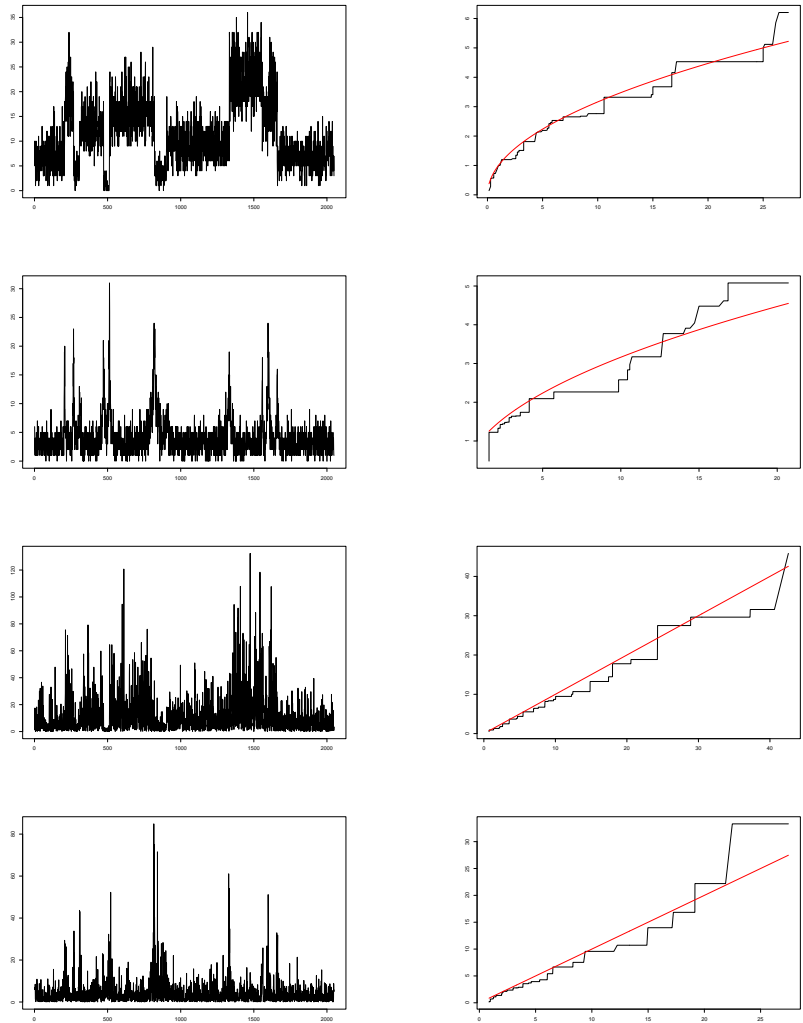

FIG 1. *Left column, from top to bottom: sample paths from the Poisson model with the* blocks *and* bumps *functions, sample paths from the exponential model with the* blocks *and* bumps *functions. Right column: corresponding estimates* $\hat{h}^{1/2}(u)$ *(step functions) and the true standard deviation functions* $h^{1/2}(u)$ *(continuous functions).*

We end this section with a brief comparison study of our estimator versus Gao's (1997) estimator for general heteroscedastic noise. The better performance

*Mean-square errors over 100 simulations for the 4 models, for Gao's method (Gao) and the Data-driven wavelet-Fisz (DdwF)*

|      | blocks exp | blocks Pois | bumps exp | bumps Pois |
| ---- | ---------- | ----------- | --------- | ---------- |
| Gao  | 7.10       | 0.93        | 3.01      | 0.91       |
| DdwF | 4.02       | 0.52        | 2.51      | 0.54       |

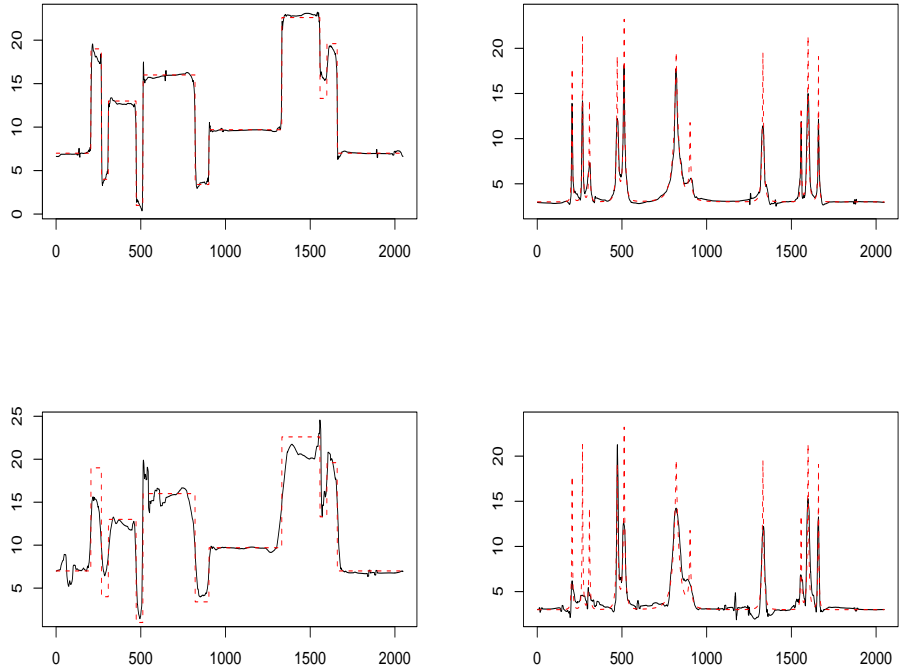

FIG 2. *Top row: estimates $\bar{\alpha}^{(h)}$ for the Poisson model (blocks, left; bumps, right). Bottom row: analogous results for the exponential model. Solid lines are the estimates, dashed lines are the true signals.*

of our method is unsurprising as Gao's method does not take into account the mean-variance relationship and thus uses less information.

## 8. Data-driven wavelet-Fisz transform

In this section, we describe a wavelet-domain variance-stabilising transform implied by our data-driven wavelet-Fisz estimation procedure.

Note that the computation of the data-driven wavelet-Fisz estimate $\bar{\alpha}^{(e)}$ can be performed in the following three steps.

1. Take a DWT of the data. For each $j = 0, \ldots, J-1$ and $k = 1, \ldots, 2^j$, divide the coefficient $Y_{j,k}$ by $\hat{h}^{1/2}(\sum_q \kappa_{j,k-q} X_q)$. Take the inverse DWT of the modified coefficients. Call the resulting vector $\tilde{X}_t$.

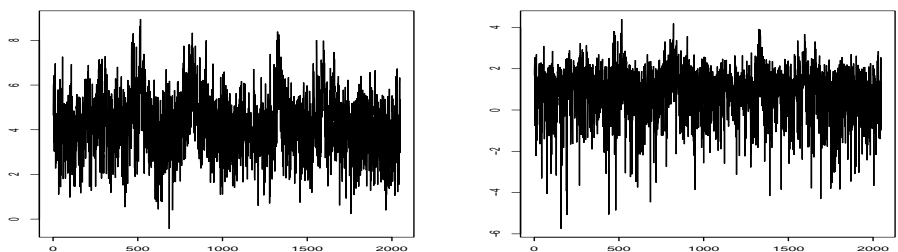

Fig 3. *Left: data-driven wavelet-Fisz transform of the exponential* bumps *dataset. Right: its log transform.*

2. Smooth $\tilde{X}_t$ by means of a standard nonlinear wavelet thresholding procedure suitable for i.i.d. Gaussian noise, using the same wavelet family as in Step 1. To be more precise, we apply the threshold $\{2\log(\#\mathcal{I}_n)\}^{1/2}$ for $(j,k) \in \mathcal{I}_n$, and set the empirical coefficients to zero for $j \in [J^*, J-1]$. Either soft or hard thresholding can be used.

3. Take the inverse transform to that described in Step 1.

We call the transform in Step 1 the *data-driven wavelet-Fisz variance-stabilising transform*. Empirically, it stabilises the variance of $X_t$ and brings its distribution closer to Gaussianity. The mechanism of the transform was already explained in the discussion underneath formula (9).

In our simulations, we found that the distribution of the "noise" in the transformed vector $\tilde{X}_t$ was the most symmetric when Haar wavelets were used. This was due to the fact that Haar wavelets are symmetric in the sense that the positive part of each Haar wavelet vector is an exact shifted version of its negative part.

Figure 3 compares the data-driven wavelet-Fisz transform (with parameters as in Section 7) of the exponential *bumps* dataset, with its logarithmic transform. Note that the latter acts as an exact variance-stabiliser, due to the multiplicative structure of the model. However, it is clear from the plot that not only does the data-driven wavelet-Fisz transform also stabilise the variance of the noise very well, but in addition it brings the distribution of the noise much closer to Gaussianity than the log transform (and does this *without knowing* the original noise distribution or the structure of the model). It also seems to bring out the shape of the underlying signal more clearly.

From a computational point of view, in view of the Gaussianising and variance-stabilising action of the data-driven wavelet-Fisz transform, the analyst wishing to find out more about the shape of the underlying signal may choose to apply *any* smoother suitable for i.i.d. Gaussian noise to the wavelet-Fisz-transformed dataset.

## 9. Discussion

We conclude with three final remarks.

*Applications.* Readers interested in the applications of the data-driven wavelet-Fisz methodology are referred to our earlier work Motakis et al. (2006) and Fryzlewicz et al. (2007), where we proposed computational procedures related to that described here and applied them to gene expression data, and solar irradiance data, respectively. Those heuristic algorithms were not accompanied by any theoretical analysis, leaving a gap which is filled by the present work.

*Density estimation.* In practice, our data-driven wavelet-Fisz estimator can also be applied to the problem of density estimation from binned data. Although this problem does not exactly fall into the class of models described by formula (2), it can be approximated, in a certain asymptotic regime, by the Poisson model, which is a sub-case of (2). We mention that Brown et al. (2007) propose a wavelet-based method for density estimation from binned data which includes a time-domain variance-stabilising transform as an initial step.

*Non-equispaced design.* Jansen et al. (2009) mention the possible use of the Haar-Fisz transform for Poisson data on graphs and in other "irregular multidimensional situations". We believe that similarly, the data-driven wavelet-Fisz methodology could also be used in such set-ups.

## Appendix A: Proof of Theorem 1

We first state an auxiliary result. Denote $\sigma_{j,k}^2 = \mathrm{Var}(\varepsilon_{j,k}) = \sum_t \psi_{j,k-t}^2 h\{\alpha(t/n)\}$ (see also formula (4)). We specify the following assumption on an arbitrary set of deterministic (non-random) thresholds $\lambda_{j,k}^{(d)}$.

**Assumption 12.**

$$\sum_{(j,k)\in\mathcal{I}_n} \left( \frac{\lambda_{j,k}^{(d)}}{\sigma_{j,k}} + 1 \right) \phi \left( \frac{\lambda_{j,k}^{(d)}}{\sigma_{j,k}} \right) = O\{n^{1/(2m+1)}\} \qquad (17)$$

$$\max_{(j,k)\in\mathcal{I}_n} \lambda_{j,k}^{(d)} = O(\log^{1/2} n), \qquad (18)$$

*where $\phi$ is the standard normal density.*

For the purposes of this section, we extend the notation $\tilde{\alpha}^{(e)}$ to mean any estimator constructed as in Section 3, using an arbitrary set of thresholds $\lambda_{j,k}$. To emphasise the dependence of $\tilde{\alpha}^{(e)}$ on $\lambda_{j,k}$, we will write $\tilde{\alpha}^{(e)}(\lambda_{j,k})$.

**Theorem 4.** *Let $\lambda_{j,k}^{(d)}$ be any non-random thresholds satisfying Assumption 12. Further, suppose that Assumptions 1, 3 and 4 hold. We have*

$$\sup_{\alpha\in\tilde{\mathcal{F}}^m} \mathbb{E}\|\tilde{\alpha}^{(e)}(\lambda_{j,k}^{(d)}) - \alpha\|_{L_2[0,1]}^2 = O\left\{(\log n/n)^{2m/(2m+1)}\right\}.$$

The proof of Theorem 4 proceeds exactly like the proof of Theorem 5.2(ii) in Neumann (1996). We omit the details.

We now define what we call "lower" and "upper" deterministic thresholds $\lambda_{j,k}^{(d,l)}$ and $\lambda_{j,k}^{(d,u)}$:

$$\lambda_{j,k}^{(d,l)} = \gamma_n h^{1/2} \left\{ \sum_q \kappa_{j,k-q} \alpha(q/n) \right\} \{2 \log(\# \mathcal{I}_n)\}^{1/2} \quad (19)$$

$$\lambda_{j,k}^{(d,u)} = C \log^{1/2} n, \quad (20)$$

where $\gamma_n$ is a sequence approaching one from below, $C$ is a generic positive constant (see Lemma 1 for the permitted range of $C$) and $\kappa_{j,\tau}$ are as in formula (6). Proving that the lower and upper thresholds satisfy Assumption 12 (and thus that Theorem 4 holds for $\tilde{\alpha}^{(e)}(\lambda_{j,k}^{(d,l)})$ and $\tilde{\alpha}^{(e)}(\lambda_{j,k}^{(d,u)})$) is a step on the way to proving Theorem 1.

**Lemma 1.** *Suppose Assumptions 1, 2 and 5 hold. If $C \geq (2\overline{h})^{1/2}$, then for all $j, k$, $\lambda_{j,k}^{(d,u)} \geq \sup_{j,k} \lambda_{j,k}^{(d,l)}$, and both $\lambda_{j,k}^{(d,l)}$ and $\lambda_{j,k}^{(d,u)}$ satisfy Assumption 12.*

PROOF. It is easy to check that if $C \geq (2\overline{h})^{1/2}$, then for all $j, k$, $\lambda_{j,k}^{(d,u)} \geq \sup_{j,k} \lambda_{j,k}^{(d,l)}$. We first check (17) for $\lambda_{j,k}^{(d,l)}$. The factor $\lambda_{j,k}^{(d,l)}/\sigma_{j,k} + 1$ only contributes a logarithmic term so we skip it. Denote $\underline{\alpha}_{j,k} = \inf\{\alpha(t/n) : t \in \mathrm{supp}(\psi_{j,k})\}$ and $\overline{\alpha}_{j,k} = \sup\{\alpha(t/n) : t \in \mathrm{supp}(\psi_{j,k})\}$. Further, let $\mathrm{TV}(f)|_A$ denote the total variation of the function $f$ measured on the set $A$. Using Assumption 2, we bound $\lambda_{j,k}^{(d,l)}/\sigma_{j,k}$ from below as follows:

$$\begin{aligned}
\frac{\lambda_{j,k}^{(d,l)}}{\sigma_{j,k}} &\geq \frac{\gamma_n h^{1/2}(\underline{\alpha}_{j,k})\{2 \log(\# \mathcal{I}_n)\}^{1/2}}{h^{1/2}(\overline{\alpha}_{j,k})} \geq \frac{\gamma_n h^{1/2}(\underline{\alpha}_{j,k})\{2 \log(\# \mathcal{I}_n)\}^{1/2}}{h^{1/2}(\underline{\alpha}_{j,k}) + H(\overline{\alpha}_{j,k} - \underline{\alpha}_{j,k})} \\
&\geq \frac{\gamma_n h^{1/2}(\underline{\alpha}_{j,k})\{2 \log(\# \mathcal{I}_n)\}^{1/2}}{h^{1/2}(\underline{\alpha}_{j,k}) + H \mathrm{TV}(\alpha)|_{\mathrm{supp}(\psi_{j,k})}} \geq \frac{\gamma_n \underline{h}^{1/2}\{2 \log(\# \mathcal{I}_n)\}^{1/2}}{\underline{h}^{1/2} + H \mathrm{TV}(\alpha)|_{\mathrm{supp}(\psi_{j,k})}},
\end{aligned}$$

where the last inequality follows from the fact that $v(x) = x/(x + a^2)$ is increasing on $[0, \infty)$. As in Neumann (1996), the proof of Lemma 6.1(ii), we have $\sum_k \mathrm{TV}(\alpha)|_{\mathrm{supp}(\psi_{j,k})} \leq O(1) \mathrm{TV}(\alpha)|_{[0,1]}$ and thus for a sequence $\omega_n \to 0$, at each scale $j$ we have

$$\#\{k : \mathrm{TV}(\alpha)|_{\mathrm{supp}(\psi_{j,k})} > \omega_n\} = O(\omega_n^{-1}). \quad (21)$$

Denote $\mathcal{D}_n = \mathcal{I}_n \cap \{(j,k) : \mathrm{TV}(\alpha)|_{\mathrm{supp}(\psi_{j,k})} > \omega_n\}$ note that by (21), at each scale $j$ at most $O(\omega_n^{-1})$ coefficients are in $\mathcal{D}_n$. Denote further $\mathcal{E}_n = \mathcal{I}_n \setminus \mathcal{D}_n$. We have

$$\sum_{\mathcal{I}_n} \phi(\lambda_{j,k}^{(d,l)}/\sigma_{j,k}) = \left( \sum_{\mathcal{D}_n} + \sum_{\mathcal{E}_n} \right) \phi(\lambda_{j,k}^{(d,l)}/\sigma_{j,k}) \leq O(\omega_n^{-1} \log n)$$

$$+ \sum_{j=1}^{J^*-1} \sum_{k=1}^{2^j} \phi\left(\frac{\gamma_n \underline{h}^{1/2}\{2\log(\#\mathcal{I}_n)\}^{1/2}}{\underline{h}^{1/2} + H\omega_n}\right) \leq O(\omega_n^{-1}\log n)$$

$$+ (2\pi)^{-1/2} \sum_{j=0}^{J^*-1} 2^{j-J^*}\gamma_n^2\left(\frac{\underline{h}^{1/2}}{\underline{h}+H\omega_n}\right)^2 = O(\omega_n^{-1}\log n)$$

$$+ O\left\{(\#\mathcal{I}_n)^{1-\gamma_n^2\left(\frac{\underline{h}^{1/2}}{\underline{h}+H\omega_n}\right)^2}\right\} = O(\omega_n^{-1}\log n) + o\{N^{1/(2m+1)}\},$$

for any $m > 0$. The last equality follows from the fact that $1 - \gamma_n^2\left(\frac{\underline{h}^{1/2}}{\underline{h}+H\omega_n}\right)^2 \to 0$. Choosing $\omega_n = \log^{-1} n$ (say), we have that (17) is satisfied irrespective of the smoothness parameter $m$. Because the thresholds $\lambda_{j,k}^{(d,u)}$ are higher than $\lambda_{j,k}^{(d,l)}$, (17) also holds for $\lambda_{j,k}^{(d,u)}$. Obviously, (18) holds for $\lambda_{j,k}^{(d,u)}$, which implies that it also holds for $\lambda_{j,k}^{(d,l)}$, since $\lambda_{j,k}^{(d,l)}$ are lower than $\lambda_{j,k}^{(d,u)}$. $\qquad\square$

We now state another auxiliary result. We first specify an assumption on an arbitrary set of *random* thresholds $\hat{\lambda}_{j,k}^{(r)}$.

**Assumption 13.**

$$\sum_{(j,k)\in\mathcal{I}_n} P(\hat{\lambda}_{j,k}^{(r)} < \lambda_{j,k}^{(d,l)}) = O(n^\nu) \qquad (22)$$

$$\sum_{(j,k)\in\mathcal{I}_n} P(\hat{\lambda}_{j,k}^{(r)} > \lambda_{j,k}^{(d,u)}) = O(n^{-2m/(2m+1)}) \qquad (23)$$

for some $\gamma_n \to 1_-$ *(see the definition of $\lambda_{j,k}^{(d,l)}$ in formula (19)), some $\nu < 1/(2m+1)$ (with $m$ given in Theorem 1), and some $C \geq (2\overline{h})^{1/2}$ (see the definition of $\lambda_{j,k}^{(d,u)}$ in formula (20)).*

**Theorem 5.** *Let $\hat{\lambda}_{j,k}^{(r)}$ be any random thresholds satisfying Assumption 13. Further, suppose that Assumptions 1, 2, 3, 4 and 5 hold. We have*

$$\sup_{\alpha\in\mathcal{F}^m} \mathbb{E}\|\tilde{\alpha}^{(e)}(\hat{\lambda}_{j,k}^{(r)}) - \alpha\|_{L_2[0,1]}^2 = O\left\{(\log n/n)^{2m/(2m+1)}\right\}.$$

The proof of Theorem 5 proceeds exactly like the proof of Theorem 6.1 in Neumann (1996). We omit the details.

In view of Theorem 5, in order to prove Theorem 1, it suffices to show that our random thresholds $\hat{\lambda}_{j,k}$, defined in formula (6), satisfy Assumption 13.

**Lemma 2.** *Suppose Assumptions 2, 3 and 5 hold. There exists a $\gamma_n \to 1_-$ and a $C \geq (2\overline{h})^{1/2}$ such that our random thresholds $\hat{\lambda}_{j,k}$, defined in formula (6), satisfy Assumption 13 for all $\nu < 1/(2m+1)$.*

PROOF. We start with (22). To shorten notation, denote $\hat{u} = \sum_q \kappa_{j,k-q} X_q$ and $u = \sum_q \kappa_{j,k-q} \alpha(q/n)$. Denote further

$$\nu_n = (1 - \gamma_n^{2\tilde{\delta}})^{1/\tilde{\delta}} \underline{h}^{\tilde{\delta}/\delta} \tilde{H}^{-1/\tilde{\delta}}.$$

Note that $\nu_n \to 0$. We have

$$P(\hat{\lambda}_{j,k} < \lambda_{j,k}^{(d,l)}) = P\{h^{1/2}(\hat{u}) < \gamma_n h^{1/2}(u)\} = P\{h^{\tilde{\delta}}(\hat{u}) < \gamma_n^{2\tilde{\delta}} h^{\tilde{\delta}}(u)\}$$

$$= P\{h^{\tilde{\delta}}(u) - h^{\tilde{\delta}}(\hat{u}) > (1 - \gamma_n^{2\tilde{\delta}}) h^{\tilde{\delta}}(u)\} \leq P(|\hat{u} - u| > \nu_n).$$

Suppose $\nu_n$ tends to zero logarithmically fast in $n$ (which is easy to ensure by placing an appropriate assumption on the speed of convergence of $\gamma_n$ to one). Then, by Lemma 8, there exists $\tilde{\epsilon} > 0$ such that

$$P(|\hat{u} - u| > \nu_n) \leq \tilde{C}_4 \exp\left(-\frac{n^{2\tilde{\epsilon}}}{2}\right).$$

Summing up over $j, k$ we obtain

$$\sum_{(j,k) \in \mathcal{I}_n} P(\hat{\lambda}_{j,k} < \lambda_{j,k}^{(d,l)}) \leq \tilde{C}_4 n \exp\left(-\frac{n^{2\tilde{\epsilon}}}{2}\right) = o(n^\nu),$$

for any $\nu$, which shows (22). The technique for showing (23) is exactly the same. We omit the details. □

The proof of Theorem 1 is complete.

## Appendix B: Proofs of results of Section 4

PROOF OF PROPOSITION 1. Let $z$ be any point such that $u = \alpha(z)$ and let $\Lambda$ be the Lipschitz constant for $\alpha$. Denote

$$B_n^z = \{t \in 1, \ldots, n : |\alpha(t/n) - \alpha(z)| \leq \Lambda |t/n - z|\}.$$

Because $\alpha(z)$ is piecewise Lipschitz-continuous of order 1, it is clear that the cardinality of $B_n^z$ is uniformly bounded from below by $cn$ where $c \in (0,1)$ and that $B_n^z$ contains those $t$ for which $t/n$ is arbitrarily close to $z$ from either the left- or the right-hand side (or both). We have

$$\sum_{t=1}^n \frac{1}{nb} K\left(\frac{\alpha(t/n) - u}{b}\right) \geq \sum_{t=1}^n \frac{1}{nb} K\left(\frac{\alpha(t/n) - \alpha(z)}{b}\right) \mathbb{I}(|\alpha(t/n) - \alpha(z)| \leq b/4)$$

$$\geq \sum_{t=1}^n \frac{1}{nb} \mathbb{I}(|\alpha(t/n) - \alpha(z)| \leq b/4) \geq \sum_{t \in B_n^z} \frac{1}{nb} \mathbb{I}(|\alpha(t/n) - \alpha(z)| \leq b/4)$$

$$\geq \sum_{t \in B_n^z} \frac{1}{nb} \mathbb{I}(\Lambda |t/n - z| \leq b/4) \geq \frac{1}{nb} nbc_1 = c_1,$$

which completes the proof. □

**Lemma 3.** *If Assumption 3 holds, then there exists a constant $\tilde{K} > 0$ such that, for $l = 3, 4, \ldots$,*

$$
\begin{aligned}
\mathbb{E}|\hat{\alpha}_t - \mathbb{E}(\hat{\alpha}_t)|^l &\leq \operatorname{Var}(\hat{\alpha}_t)(l!)^{2+\gamma} \left\{ \frac{\tilde{K}}{(2M+1)^{1/2}} \right\}^{l-2} \\
&= (2M+1)^{-2} \sum_{p=t-M}^{t+M} h\{\alpha(p/n)\}(l!)^{2+\gamma} \left\{ \frac{\tilde{K}}{(2M+1)^{1/2}} \right\}^{l-2}.
\end{aligned}
$$

PROOF. Let $C$ be a generic positive constant. Recall that by Stirling's approximation,

$$
l^l \leq \exp\{l - 1/(12l+1)\} l! (2\pi l)^{-1/2} \tag{24}
$$

Using the Rosenthal inequality (Rosenthal (1970); Johnson et al. (1985)), and then (24) and Assumption 3, we have

$$
\mathbb{E}|\hat{\alpha}_t - \mathbb{E}(\hat{\alpha}_t)|^l = (2M+1)^{-l} \mathbb{E}\left| \sum_{p=t-M}^{t+M} X_p - \alpha(p/n) \right|^l \leq (2M+1)^{-l} \frac{C^l l^l}{\log^l l}
$$

$$
\times \max\left\{ \sum_{p=t-M}^{t+M} \mathbb{E}|X_p - \alpha(p/n)|^l, \left( \sum_{p=t-M}^{t+M} h\{\alpha(p/n)\} \right)^{l/2} \right\}
$$

$$
\leq (2M+1)^{-l} l! \, p(l)
$$

$$
\times \max\left\{ (l!)^{1+\gamma} K^{l-2} \sum_{p=t-M}^{t+M} h\{\alpha(p/n)\}, \left( \sum_{p=t-M}^{t+M} h\{\alpha(p/n)\} \right)^{l/2} \right\} \tag{25}
$$

where

$$
p(l) = C^l \exp\{l - 1/(12l+1)\}(2\pi l)^{-1/2} \log^{-l} l.
$$

Noting that $p(l) \leq \operatorname{const}^l$, we observe that (25) can comfortably be bounded by

$$
(2M+1)^{-2} \sum_{p=t-M}^{t+M} h\{\alpha(p/n)\}(l!)^{2+\gamma} \left\{ \frac{\tilde{K}}{(2M+1)^{1/2}} \right\}^{l-2}
$$

$$
= \operatorname{Var}(\hat{\alpha}_t)(l!)^{2+\gamma} \left\{ \frac{\tilde{K}}{(2M+1)^{1/2}} \right\}^{l-2}
$$

for a constant $\tilde{K} > 0$, which completes the proof. □

**Lemma 4.** *For a nonnegative random variable $X$ and $l \geq 1$, we have*

$$
\mathbb{E}|X - \mathbb{E}(X)|^l \leq 2\mathbb{E}(X^l).
$$

PROOF. Noting that for $c \geq 0$, we have $|x - c|^l \leq |x|^l + c^l$, we obtain

$$\mathbb{E}|X - \mathbb{E}(X)|^l \leq \mathbb{E}(X^l) + \{\mathbb{E}(X)\}^l \leq 2\mathbb{E}(X^l),$$

where the last step uses Jensen's inequality.                                  □

**Lemma 5.** *Denote* $Z_t = |\hat{\alpha}_t - \mathbb{E}(\hat{\alpha}_t)|$*. Under Assumptions 3 and 8, we have*

$$\mathbb{E}|Z_t - \mathbb{E}(Z_t)|^l \leq \text{Var}(Z_t)(l!)^{2+\gamma} \left\{ \frac{\bar{K}}{(2M+1)^{1/2}} \right\}^{l-2}$$

*for* $l = 3, 4, \ldots$*, where* $\bar{K}$ *is a positive constant.*

PROOF. Using Lemma 4 with $X = Z_t$ and then Lemma 3 and Assumption 8, we have

$$\begin{aligned}
\mathbb{E}|Z_t - \mathbb{E}(Z_t)|^l &\leq& 2\mathbb{E}|\hat{\alpha}_t - \mathbb{E}(\hat{\alpha}_t)|^l \leq \text{Var}(\hat{\alpha}_t)(l!)^{2+\gamma} \left\{ \frac{2\tilde{K}}{(2M+1)^{1/2}} \right\}^{l-2} \\
&\leq& \text{Var}(Z_t)(l!)^{2+\gamma} \left\{ \frac{\bar{K}}{(2M+1)^{1/2}} \right\}^{l-2}
\end{aligned}$$

for $\bar{K} = 2C_2\tilde{K}$.                                               □

**Lemma 6.** *Denote* $\tilde{\varepsilon}_t = \varepsilon_t^2 - h\{\alpha(t/n)\}$*. Under Assumptions 2, 3 and 7, we have*

$$\mathbb{E}|\tilde{\varepsilon}_t|^l \leq (l!)^{3+3\gamma} K_1^{l-2} \text{Var}(\tilde{\varepsilon}_t),$$

*for* $l = 3, 4, \ldots$*, where*

$$\text{Var}(\tilde{\varepsilon}_t) = \text{Var}(\varepsilon_t^2) \leq 24^{1+\gamma} K^2 h\{\alpha(t/n)\} - h^2\{\alpha(t/n)\}.$$

PROOF. Applying Lemma 4 to $\varepsilon_t^2$ and then using Assumption 3, we get

$$\mathbb{E}|\tilde{\varepsilon}_t|^l \leq 2\mathbb{E}|\varepsilon_t|^{2l} \leq 2\{(2l)!\}^{1+\gamma} K^{2l-2} h\{\alpha(t/n)\}.$$

Using the fact that $(2l)! \leq 4(l!)^3$ and Assumption 7, we bound the above by $(l!)^{3+3\gamma} K_1^{l-2} \text{Var}(\tilde{\varepsilon}_t)$, which completes the proof.               □

**Lemma 7.** *Suppose that Assumptions 1, 2, 3 and 8 hold, and that* $\delta = \delta_n$ *is such that*

$$\delta n^{-1/2} - Mn^{-1/2} - n^{1/2}M^{-1/2} \to \infty.$$

*Let* $b_n$ *be any fixed sequence such that* $b_n = o\{(n/M)^{1/(6+4\gamma)}\}$*. In the asymptotic limit, as* $n, M, n/M \to \infty$*, we have*

$$P\left( \sum_{t=1}^{n} |\hat{\alpha}_t - \alpha(t/n)| \geq \delta \right) \leq C_3(2M+1)\{1 - \Phi(\min(a_n, b_n))\}, \qquad (26)$$

*where* $C_3$ *is a positive constant,* $\Phi$ *is the cdf of the standard normal, and the sequence* $a_n$ *satisfies*

$$a_n = O(\delta n^{-1/2} - Mn^{-1/2} - n^{1/2}M^{-1/2}).$$

PROOF. We first note that if Assumption 1 (i) holds, then

$$\frac{1}{2M+1} \sum_{t=i,i+1+2M,i+2+4M,\dots} \sum_{p=t-M}^{t+M} |\alpha(p/n) - \alpha(t/n)| \le C_1 \qquad (27)$$

uniformly over $i$ and $M$, where $C_1$ is a positive constant. Summing up both sides of equation (27) over $i$, we have $\sum_{t=1}^{n} |\alpha(t/n) - \mathbb{E}(\hat{\alpha}_t)| \le C_1(2M+1)$. We bound

$$
\begin{aligned}
P\left( \sum_{t=1}^{n} |\hat{\alpha}_t - \alpha(t/n)| \ge \delta \right) &\le P\left( \sum_{t=1}^{n} |\hat{\alpha}_t - \mathbb{E}(\hat{\alpha}_t)| + |\mathbb{E}(\hat{\alpha}_t) - \alpha(t/n)| \ge \delta \right) \\
&\le P\left( \sum_{t=1}^{n} |\hat{\alpha}_t - \mathbb{E}(\hat{\alpha}_t)| \ge \delta - C_1(2M+1) \right) \quad (28)
\end{aligned}
$$

Denote $\delta' := \delta - C_1(2M+1)$. The sequence $\{\hat{\alpha}_t\}_t$ is $(2M+1)-$dependent. To avoid complications which arise in deriving exponential inequalities for dependent sequences, we split it into independent sequences as follows. Rewriting the LHS of (28) as

$$P\left( \sum_i \sum_{t=i,i+1+2M,i+2+4M,\dots} |\hat{\alpha}_t - \mathbb{E}(\hat{\alpha}_t)| \ge \delta' \right), \qquad (29)$$

and using the fact that $a_1 + \dots + a_m \ge \delta \implies \exists i \quad a_i \ge \delta/m$, as well as the Bonferroni inequality, we bound (29) by

$$(2M+1) \max_i P\left( \sum_{t=i,i+1+2M,i+2+4M,\dots} |\hat{\alpha}_t - \mathbb{E}(\hat{\alpha}_t)| \ge \delta'/(2M+1) \right).$$

We drop the range of $t$ to shorten notation, and assume without loss of generality that there are exactly $n/(2M+1)$ terms in the above sum. Denoting $Z_t = |\hat{\alpha}_t - \mathbb{E}(\hat{\alpha}_t)|$, we bound the above by

$$(2M+1) \max_i P\left( \sum_t Z_t - \mathbb{E}(Z_t) \ge \delta'/(2M+1) - \sum_t \mathbb{E}(Z_t) \right). \qquad (30)$$

We first assess $\sum_t \mathbb{E}(Z_t)$:

$$\sum_t \mathbb{E}(Z_t) \le \sum_t \{\mathrm{var}(\hat{\alpha}_t)\}^{1/2} \le n\overline{h}^{1/2}(2M+1)^{-3/2}.$$

Denote

$$\delta'' = \frac{\delta'/(2M+1) - n\overline{h}^{1/2}(2M+1)^{-3/2}}{\{\sum_t \mathrm{Var}(Z_t)\}^{1/2}}$$

$$\xi = \frac{\sum_t Z_t - \mathbb{E}(Z_t)}{\{\sum_t \text{Var}(Z_t)\}^{1/2}}.$$

Note that by Assumption 8, we have $\delta'' = O(\delta n^{-1/2} - Mn^{-1/2} - n^{1/2}M^{-1/2})$ and, by assumptions of the lemma, $\delta'' \to \infty$. We now apply Theorem 1 and its Corollary from Rudzkis et al. (1978) to the standardised sum $\xi$. By Lemma 5, in our case the quantity $\Delta_n$ from the above Theorem takes the form

$$\Delta_n = \frac{\{\sum_t \text{Var}(Z_t)\}^{1/2}}{2\max\{\bar{K}(2M+1)^{-1/2}, \max_t\{\text{Var}(Z_t)\}^{1/2}\}},$$

which, by Assumption 8, is of order $O\{(n/M)^{1/2}\}$. Recall that $b_n = o\{(n/M)^{1/(6+4\gamma)}\}$. Using the above Theorem, we bound (30) by

$$
\begin{aligned}
(2M+1)\max_i P(\xi \geq \delta'') &\leq (2M+1)\max_i P\{\xi \geq \min(\delta'', b_n)\} \\
&\leq C_3(2M+1)\{1 - \Phi(\min(\delta'', b_n))\},
\end{aligned}
$$

which completes the proof. $\qquad\square$

**Lemma 8.** *Let the constants $\kappa_{j,\tau}$ satisfy Assumption 5. Suppose Assumptions 2 and 3 hold. Recall that $n = 2^J$. Let $b_n$ be any fixed sequence such that $b_n = o\left(2^{\frac{J-j}{2(1+\max(\gamma,1))}}\right)$ uniformly over $0 \leq j \leq J^* - 1$. In the asymptotic limit, as $n \to \infty$, and uniformly over $0 \leq j \leq J^* - 1$, we have*

$$P\left(\left|\sum_\tau \kappa_{j,\tau}\{X_\tau - \alpha(\tau/n)\}\right| > \delta\right) \leq C_4(1 - \Phi(\min(a_n, b_n))),$$

*where $C_4$ is a positive constant independent of $j$, $\Phi$ is the cdf of the standard normal, and the sequence $a_n$ satisfies*

$$a_n = O(\delta \min_j 2^{(J-j)/2}).$$

PROOF. Denote $\tilde{X}_\tau = \kappa_{j,\tau}\{X_\tau - \alpha(\tau/n)\}$. We have $\text{Var}(\tilde{X}_\tau) = \kappa_{j,\tau}^2 h\{\alpha(\tau/n)\}$ and, by Assumption 3,

$$\mathbb{E}|\tilde{X}_\tau|^l \leq (l!)^{1+\gamma}(K\kappa_{j,\tau})^{l-2}\text{Var}(\tilde{X}_\tau) \leq (l!)^{1+\gamma}(K\max_k \kappa_{j,k})^{l-2}\text{Var}(\tilde{X}_\tau) \quad (31)$$

for $l = 3, 4, \ldots$. Denote $\xi = \sum_\tau \tilde{X}_\tau / \{\sum_\tau \text{Var}(\tilde{X}_\tau)\}^{1/2}$ and $\delta' = \delta/\{\sum_\tau \text{Var}(\tilde{X}_\tau)\}^{1/2}$. By Assumption 5, $\delta' = O(\delta 2^{(J-j)/2})$. We now apply Theorem 1 and its Corollary from Rudzkis et al. (1978) to the standardised sum $\xi$. By (31), in our case the quantity $\Delta_n$ from the above Theorem takes the form

$$\Delta_n = \Delta_{n,j} = \frac{\{\sum_\tau \text{Var}(\tilde{X}_\tau)\}^{1/2}}{2\max\{K\max_k \kappa_{j,k}, \max_\tau\{\text{Var}(\tilde{X}_\tau)\}^{1/2}\}},$$

which, by Assumption 5, is of order $O(2^{(J-j)/2})$. Recall that $b_n = o\left(2^{\frac{J-j}{2(1+\max(\gamma,1))}}\right)$. Using the above Theorem, we bound

$$P(|\xi| > \delta') \leq P(|\xi| > \min(\delta', b_n)) \leq C_4(1 - \Phi(\min(\delta', b_n))),$$

which completes the proof. □

We define $\omega_{nt}(u) = W_{nt}(u) / \sum_{t=1}^n W_{nt}(u)$. Obviously $\sum_{t=1}^n \omega_{nt}(u) = 1$, and if Assumptions 6 and 10 hold, then $\omega_{nt}(u) \leq \dot{K}/(c_1 nb) = O(n^{-1}b^{-1})$ so that $\sum_{t=1}^n \omega_{nt}^2(u) \leq O(n^{-1}b^{-2})$. By Cauchy inequality, we also have $1 = \{\sum_{t=1}^n \omega_{nt}(u)\}^2 \leq n \sum_{t=1}^n \omega_{nt}^2(u)$, which implies $\sum_{t=1}^n \omega_{nt}^2(u) \geq 1/n$. Summarising the above bounds,

$$O(n^{-1}b^{-2}) \geq \sum_{t=1}^n \omega_{nt}^2(u) \geq 1/n. \qquad (32)$$

**Lemma 9.** *Suppose Assumptions 2, 3, 6, 7 and 10 hold. Let $b_n$ be any fixed sequence such that $b_n = o\left((nb^2)^{1/(10+12\gamma)}\right)$. In the asymptotic limit, as $nb \to \infty$, $b \to 0$, we have*

$$P\left(\left|\sum_{t=1}^n \omega_{nt}(u)(\varepsilon_t^2 - h\{\alpha(t/n)\})\right| \geq \delta\right) \leq C_5(1 - \Phi(\min(a_n, b_n))),$$

*where $C_5$ is a positive constant, $\Phi$ is the cdf of the standard normal, and the sequence $a_n$ satisfies*

$$a_n = O(\delta n^{1/2}b).$$

PROOF. Denote $\bar{\varepsilon}_t = \omega_{nt}(u)(\varepsilon_t^2 - h\{\alpha(t/n)\})$. By Lemma 6 and Assumption 10,

$$\mathbb{E}|\bar{\varepsilon}_t|^l \leq \omega_{nt}^l(u)(l!)^{3+3\gamma} K_1^{l-2} \text{Var}(\varepsilon_t^2 - h\{\alpha(t/n)\}) \leq (l!)^{3+3\gamma} (K')^{l-2} \text{Var}(\bar{\varepsilon}_t),$$

for $l = 3, 4, \ldots$, where $K' = K_1 \dot{K}/(c_1 nb) = O(n^{-1}b^{-1})$. Denote $\xi = \sum_t \bar{\varepsilon}_t / \{\sum_t \text{Var}(\bar{\varepsilon}_t)\}^{1/2}$ and $\delta' = \delta / \{\sum_t \text{Var}(\bar{\varepsilon}_t)\}^{1/2}$. By (32), $\delta' \geq O(\delta n^{1/2}b)$. We now apply Theorem 1 and its Corollary from Rudzkis et al. (1978) to the standardised sum $\xi$. In our case the quantity $\Delta_n$ from the above Theorem takes the form

$$\Delta_n = \frac{\{\sum_t \text{Var}(\bar{\varepsilon}_t)\}^{1/2}}{2\max\{K', \max_t\{\text{Var}(\bar{\varepsilon}_t)\}^{1/2}\}},$$

which, by Assumption 10 and (32) is of order at least $O(n^{1/2}b)$. Recall that $b_n = o\left((nb^2)^{1/(10+12\gamma)}\right)$. Using the above Theorem, we bound

$$P(|\xi| \geq \delta') \leq P(|\xi| \geq \min(\delta', b_n)) \leq C_5(1 - \Phi(\min(\delta', b_n))),$$

which completes the proof. □

**Lemma 10.** *Suppose that Assumptions 1, 2, 3, 6, 8 hold, and that the constants $\kappa_{j,\tau}$ satisfy Assumption 5. Let $b_n$, $d_n$ be any fixed sequences s.t. $b_n =$* $o\{(n/M)^{1/(6+4\gamma)}\}$ *and* $d_n = o\left(\min_j 2^{\frac{J-j}{2(1+\max(\gamma,1))}}\right)$, *where* $0 \leq j \leq J^* - 1$. *Let* $\delta$ *be such that*

$$\delta n^{1/2} b^2 - M n^{-1/2} - n^{1/2} M^{-1/2} \to \infty,$$

*as* $M, n, n/M \to \infty$ *and* $b \to 0$. *In the asymptotic limit, as* $n, M, n/M \to \infty$ *and* $b \to 0$, *uniformly over* $j, k$, *we have*

$$P\left\{\sum_{t=1}^{n} W_{nt}\left(\sum_q \kappa_{j,k-q}\alpha(q/n)\right) - \hat{W}_{nt}\left(\sum_q \kappa_{j,k-q}X_q\right) \geq \delta\right\}$$

$$\leq C_3(2M+1)\{1 - \Phi(\min(a_n, b_n))\} + C_4\{1 - \Phi(\min(c_n, d_n))\}, \quad (33)$$

*where* $C_3$, $C_4$ *are as in Lemmas* 7 *and* 8 *and* $\Phi$ *is the cdf of standard normal.*

$$\begin{aligned} a_n &= O(\delta n^{1/2} b^2 - M n^{-1/2} - n^{1/2} M^{-1/2}) \\ c_n &= O(\delta b^2 \min_j 2^{(J-j)/2}). \end{aligned}$$

PROOF. Using the Lipschitz-continuity of the kernel function $K(\cdot)$, we have

$$\left| W_{nt}\left(\sum_q \kappa_{j,k-q}\alpha(q/n)\right) - \hat{W}_{nt}\left(\sum_q \kappa_{j,k-q}X_q\right) \right|$$

$$\leq \frac{L}{nb^2}\left\{|\alpha(t/n) - \hat{\alpha}_t| - \left|\sum_q \kappa_{j,k-q}\{X_q - \alpha(q/n)\}\right|\right\}$$

Thus, we bound the probability on the LHS of (33) by

$$P\left(\sum_{t=1}^{n}|\alpha(t/n) - \hat{\alpha}_t| \geq \frac{\delta n b^2}{2L}\right) + P\left(\left|\sum_q \kappa_{j,k-q}\{X_q - \alpha(q/n)\}\right| \geq \frac{\delta b^2}{2L}\right).$$

Lemmas 7 and 8 yield the result. □

PROOF OF THEOREM 2. To shorten notation, denote $u = \sum_q \kappa_{j,k-q}\alpha(q/n)$ and $\hat{u} = \sum_q \kappa_{j,k-q}X_q$. By Assumption 11, $u \in \text{range}\{\alpha(z)\}$. We bound

$$P(|\hat{h}(\hat{u}) - h(u)| \geq \delta) \leq P(|\hat{h}(\hat{u}) - \tilde{h}(u)| \geq \delta/2) + P(|\tilde{h}(u) - h(u)| \geq \delta/2) =: I + II,$$

where

$$\tilde{h}(u) = \frac{\sum_{t=1}^{n} W_{nt}(u)\varepsilon_t^2}{\sum_{t=1}^{n} W_{nt}(u)}.$$

We first consider $I$. Again to shorten notation, denote

$$A = \sum_{t=1}^{n} \hat{W}_{nt}(\hat{u})\hat{\varepsilon}_t^2 \sum_{t=1}^{n} W_{nt}(u) \qquad\qquad B = \sum_{t=1}^{n} W_{nt}(u)\varepsilon_t^2 \sum_{t=1}^{n} \hat{W}_{nt}(\hat{u})$$

$$C = \sum_{t=1}^{n} W_{nt}(u)\varepsilon_t^2 \sum_{t=1}^{n} W_{nt}(u) \qquad\qquad D = \sum_{t=1}^{n} \hat{W}_{nt}(\hat{u}) \sum_{t=1}^{n} W_{nt}(u)$$

and also

$$E = \left\{ \left| \sum_{t=1}^{n} \hat{W}_{nt}(\hat{u}) - \sum_{t=1}^{n} W_{nt}(u) \right| \geq \delta_1 \right\}.$$

We bound

$$
\begin{aligned}
I &= P(|A - B| \geq \delta D/2) \\
&= P(|A - B| \geq \delta D/2 | E)P(E) + P(|A - B| \geq \delta D/2 | E^c)P(E^c) \\
&\leq P(E) + P\left( |A - B| \geq \frac{\delta}{2} \sum_{t=1}^{n} W_{nt}(u) \left\{ \sum_{t=1}^{n} W_{nt}(u) - \delta_1 \right\} \middle| E^c \right) P(E^c) \\
&\leq P(E) + P\left( |A - B| \geq \frac{\delta}{2} \sum_{t=1}^{n} W_{nt}(u) \left\{ \sum_{t=1}^{n} W_{nt}(u) - \delta_1 \right\} \right) \\
&\leq P(E) + P\left( |A - B| \geq \delta' \sum_{t=1}^{n} W_{nt}(u) \right).
\end{aligned}
$$

By Lemma 10,

$$P(E) \leq C_3(2M + 1)\left\{1 - \Phi(\min(a_n, b_n))\right\} + C_4\left\{1 - \Phi(\min(c_n, d_n))\right\}. \quad (34)$$

We bound

$$
\begin{aligned}
&P\left( |A - B| \geq \delta' \sum_{t=1}^{n} W_{nt}(u) \right) \\
&\leq P\left( |A - C| \geq \frac{\delta'}{2} \sum_{t=1}^{n} W_{nt}(u) \right) + P\left( |C - B| \geq \frac{\delta'}{2} \sum_{t=1}^{n} W_{nt}(u) \right) =: I_1 + I_2.
\end{aligned}
$$

Turning first to $I_1$, we have

$$
\begin{aligned}
I_1 &= P\left( \left| \sum_{t=1}^{n} \hat{W}_{nt}(\hat{u})\hat{\varepsilon}_t^2 - \sum_{t=1}^{n} W_{nt}(u)\varepsilon_t^2 \right| \geq \frac{\delta'}{2} \right) \\
&\leq P\left( \left| \sum_{t=1}^{n} \varepsilon_t^2 \left\{ W_{nt}(u) - \hat{W}_{nt}(\hat{u}) \right\} \right| \geq \frac{\delta'}{4} \right) \\
&\quad + P\left( \left| \sum_{t=1}^{n} \hat{W}_{nt}(\hat{u}) \left\{ \varepsilon_t^2 - \hat{\varepsilon}_t^2 \right\} \right| \geq \frac{\delta'}{4} \right) =: I_{11} + I_{12}.
\end{aligned}
$$

We first consider $I_{11}$. Denote

$$F = \{ \forall t = 1, \ldots, n \quad |\varepsilon_t| \leq a \log^d n \},$$

where $a$, $d$ are constants from Assumption 9. Using Assumption 9, the assumptions of this theorem, and Lemma 10, we bound

$$
I_{11} \leq P(F^c) + P\left( \left| \sum_{t=1}^{n} a^2 \log^{2d} n \left\{ W_{nt}(u) - \hat{W}_{nt}(\hat{u}) \right\} \right| \geq \frac{\delta'}{4} \right)
$$

$$
\begin{aligned}
&= & O(n^{-2}) + P\left(\left|\sum_{t=1}^{n} W_{nt}(u) - \hat{W}_{nt}(\hat{u})\right| \geq \frac{\delta'}{4a^2\log^{2d} n}\right) \\
&\leq & O(n^{-2}) + C_3(2M+1)\{1 - \Phi(\min(e_n, b_n))\} \\
&+ & C_4\{1 - \Phi(\min(f_n, d_n))\}.
\end{aligned}
\tag{35}
$$

We now consider $I_{12}$. Denote $p(\log^d n) = 3a\log^d n + \overline{\alpha}$. Noting that $\varepsilon_t - \hat{\varepsilon}_t = \hat{\alpha}_t - \alpha(t/n)$, we obtain

$$
\begin{aligned}
|I_{12}| &\leq P\left(\frac{\dot{K}}{nb}\sum_{t=1}^{n}|\varepsilon_t + \hat{\varepsilon}_t|\,|\hat{\alpha}_t - \alpha(t/n)| \geq \frac{\delta'}{4}\right) \\
&\leq P\left(\frac{\dot{K}}{nb}\sum_{t=1}^{n}\left\{2|\varepsilon_t| + \frac{1}{2M+1}\sum_{q=t-M}^{t+M}|\varepsilon_q|\right.\right. \\
&\quad + \left.\left.\left|\alpha(t/n) - \frac{1}{2M+1}\sum_{q=t-M}^{t+M}\alpha(q/n)\right|\right\}|\hat{\alpha}_t - \alpha(t/n)| \geq \frac{\delta'}{4}\right) \\
&\leq P(F^c) + P\left(\sum_{t=1}^{n}|\hat{\alpha}_t - \alpha(t/n)| \geq \frac{\delta' nb}{4\dot{K}p(\log^d n)}\right).
\end{aligned}
$$

Using Assumption 9, the assumptions of this theorem, and Lemma 7, we bound the above by

$$
|I_{12}| \leq O(n^{-2}) + C_3(2M+1)\{1 - \Phi(\min(g_n, b_n))\}.
\tag{36}
$$

We now consider $I_2$. We have

$$
|I_2| \leq P(F^c) + P\left(\left|\sum_{t=1}^{n} W_{nt}(u) - \hat{W}_{nt}(\hat{u})\right| > \frac{\delta'}{2a^2\log^{2d} n}\right)
$$

Using Assumption 9, the assumptions of this theorem, and Lemma 10, the bound is the same as that for $I_{11}$:

$$
|I_2| \leq O(n^{-2}) + C_3(2M+1)\{1 - \Phi(\min(e_n, b_n))\} + C_4\{1 - \Phi(\min(f_n, d_n))\}.
\tag{37}
$$

We finally turn to $II$. We define

$$
\bar{h}(u) := \mathbb{E}(\tilde{h}(u)) = \frac{\sum_{t=1}^{n} W_{nt}(u)h\{\alpha(t/n)\}}{\sum_{t=1}^{n} W_{nt}(u)}.
$$

Note that

$$
\begin{aligned}
|\bar{h}(u) - h(u)| &\leq & \frac{\sum_{t=1}^{n} W_{nt}(u)|h\{\alpha(t/n)\} - h(u)|}{\sum_{t=1}^{n} W_{nt}(u)} \\
&\leq & \frac{H\sum_{t=1}^{n} W_{nt}(u)|\alpha(t/n) - u|}{\sum_{t=1}^{n} W_{nt}(u)} \leq \frac{bH}{2},
\end{aligned}
$$

where the last inequality comes from the fact that $W_{nt}(u)$ is supported on $[\alpha(t/n) - b/2, \alpha(t/n) + b/2]$. We bound

$$|II| \leq P(|\tilde{h}(u) - \bar{h}(u)| + |\bar{h}(u) - h(u)| \geq \delta/2) \leq P(|\tilde{h}(u) - \bar{h}(u)| \geq \delta/2 - bH/2),$$

which, by the assumptions of this theorem and by Lemma 9 is bounded by

$$|II| \leq C_5(1 - \Phi(\min(h_n, i_n))). \tag{38}$$

Combining (34), (35), (36), (37) and (38) yields the result. □

PROOF OF COROLLARY 1. Denoting by $\tilde{\phi}(x)$ the standard normal pdf, and recalling that for large $x$, we have $1 - \Phi(x) \leq \tilde{\phi}(x)$, Corollary 1 is a direct consequence of Theorem 2 and the discussion directly underneath it. □

## Appendix C: Proof of Theorem 3

In view of Theorem 5, it suffices to show that our thresholds $\tilde{\lambda}_{j,k}$ satisfy Assumption 13. The following lemma holds.

**Lemma 11.** *Suppose that Assumptions 1, 2, 3, 5, 6, 7, 8, 9, 10 and 11 hold. There exists a $\gamma_n \to 1_-$ and a $C \geq (2\tilde{h})^{1/2}$ such that our random thresholds $\tilde{\lambda}_{j,k}$, defined in formula (16), satisfy Assumption 13 for all $-1 \leq \nu < 1/(2m + 1)$.*

PROOF. We start with (22). To shorten notation, denote $\hat{u} = \sum_q \kappa_{j,k-q} X_q$ and $u = \sum_q \kappa_{j,k-q} \alpha(q/n)$. We have

$$P(\tilde{\lambda}_{j,k} < \lambda_{j,k}^{(d,l)}) = P\{\hat{h}^{1/2}(\hat{u}) < \gamma_n h^{1/2}(u)\} = P\{\hat{h}(\hat{u}) < \gamma_n^2 h(u)\} =$$

$$P\{h(u) - \hat{h}(\hat{u}) > (1 - \gamma_n^2)h(u)\} \leq P\{|h(u) - \hat{h}(\hat{u})| > (1 - \gamma_n^2)\underline{h}\}.$$

Suppose $\gamma_n^2$ converges to one logarithmically fast in $n$. Then, by Corollary 1, the above probability can uniformly be bounded by $O(n^{-2})$. Summing over $j, k$, we obtain

$$\sum_{(j,k) \in \mathcal{I}_n} P(\tilde{\lambda}_{j,k} < \lambda_{j,k}^{(d,l)}) \leq O(n^{-1}) \leq O(n^\nu),$$

which shows (22). The technique for showing (23) is exactly the same. We omit the details. □

## Acknowledgements

I wish to thank Véronique Delouille and Jean-François Hochedez for their hospitality and financial support during my visit to the Royal Observatory of Belgium (Brussels), where this work was initiated. I am also grateful to the Associate Editor and the Referee for their comments which led to an improved version of this paper.